
%
\magnification1200
\pretolerance=100
\tolerance=200
\hbadness=1000
\vbadness=1000
\linepenalty=10
\hyphenpenalty=50
\exhyphenpenalty=50
\binoppenalty=700
\relpenalty=500
\clubpenalty=5000
\widowpenalty=5000
\displaywidowpenalty=50
\brokenpenalty=100
\predisplaypenalty=7000
\postdisplaypenalty=0
\interlinepenalty=10
\doublehyphendemerits=10000
\finalhyphendemerits=10000
\adjdemerits=160000
\uchyph=1
\delimiterfactor=901
\hfuzz=0.1pt
\vfuzz=0.1pt
\overfullrule=5pt
\hsize=146 true mm
\vsize=8.9 true in
\maxdepth=4pt
\delimitershortfall=.5pt
\nulldelimiterspace=1.2pt
\scriptspace=.5pt
\normallineskiplimit=.5pt
\mathsurround=0pt
\parindent=20pt
\catcode`\_=11
\catcode`\_=8
\normalbaselineskip=12pt
\normallineskip=1pt plus .5 pt minus .5 pt
\parskip=6pt plus 3pt minus 3pt
\abovedisplayskip = 12pt plus 5pt minus 5pt
\abovedisplayshortskip = 1pt plus 4pt
\belowdisplayskip = 12pt plus 5pt minus 5pt
\belowdisplayshortskip = 7pt plus 5pt
\normalbaselines
\smallskipamount=\parskip
 \medskipamount=2\parskip
 \bigskipamount=3\parskip
\jot=3pt
%
%
\def\ref#1{\par\noindent\hangindent2\parindent
 \hbox to 2\parindent{#1\hfil}\ignorespaces}
%
%
\font\typd=cmbx10 scaled \magstep2   
\font\typf=cmcsc10                   
\font\tenss=cmss10
\font\sevenss=cmss8 at 7pt
\font\fivess=cmss8 at 5pt
\newfam\ssfam %
\textfont\ssfam=\tenss
\scriptfont\ssfam=\sevenss
\scriptscriptfont\ssfam=\fivess
%
%
%
%
%
%
%
%
%
\catcode`\_=11
\def\suf_fix{}
\def\scaled_rm_box#1{%
 \relax
 \ifmmode
   \mathchoice
    {\hbox{\tenrm #1}}%
    {\hbox{\tenrm #1}}%
    {\hbox{\sevenrm #1}}%
    {\hbox{\fiverm #1}}%
 \else
  \hbox{\tenrm #1}%
 \fi}
\def\suf_fix_def#1#2{\expandafter\def\csname#1\suf_fix\endcsname{#2}}
\def\I_Buchstabe#1#2#3{%
 \suf_fix_def{#1}{\scaled_rm_box{I\hskip-0.#2#3em #1}}
}
\def\rule_Buchstabe#1#2#3#4{%
 \suf_fix_def{#1}{%
  \scaled_rm_box{%
   \hbox{%
    #1%
    \hskip-0.#2em%
    \lower-0.#3ex\hbox{\vrule height1.#4ex width0.07em }%
   }%
   \hskip0.50em%
  }%
 }%
}
\I_Buchstabe B22
\rule_Buchstabe C51{34}
\I_Buchstabe D22
\I_Buchstabe E22
\I_Buchstabe F22
\rule_Buchstabe G{525}{081}4
\I_Buchstabe H22
\I_Buchstabe I20
\I_Buchstabe K22
\I_Buchstabe L20
\I_Buchstabe M{20em }{I\hskip-0.35}
\I_Buchstabe N{20em }{I\hskip-0.35}
\rule_Buchstabe O{525}{095}{45}
\I_Buchstabe P20
\rule_Buchstabe Q{525}{097}{47}
\I_Buchstabe R21 
\rule_Buchstabe U{45}{02}{54}
\suf_fix_def{Z}{\scaled_rm_box{Z\hskip-0.38em Z}}
\catcode`\"=12
\newcount\math_char_code
\def\suf_fix_math_chars_def#1{%
 \ifcat#1A
  \expandafter\math_char_code\expandafter=\suf_fix_fam
  \multiply\math_char_code by 256
  \advance\math_char_code by `#1
  \expandafter\mathchardef\csname#1\suf_fix\endcsname=\math_char_code
  \let\next=\suf_fix_math_chars_def
 \else
  \let\next=\relax
 \fi
 \next}
%
%
%
%
\def\font_fam_suf_fix#1#2 #3 {%
 \def\suf_fix{#2}
 \def\suf_fix_fam{#1}
 \suf_fix_math_chars_def #3.
}
\font_fam_suf_fix
 0rm
 ABCDEFGHIJKLMNOPQRSTUVWXYZabcdefghijklmnopqrstuvwxyz
\font_fam_suf_fix
 2scr
 ABCDEFGHIJKLMNOPQRSTUVWXYZ
\font_fam_suf_fix
 \slfam sl
 ABCDEFGHIJKLMNOPQRSTUVWXYZabcdefghijklmnopqrstuvwxyz
\font_fam_suf_fix
 \bffam bf
 ABCDEFGHIJKLMNOPQRSTUVWXYZabcdefghijklmnopqrstuvwxyz
\font_fam_suf_fix
 \ttfam tt
 ABCDEFGHIJKLMNOPQRSTUVWXYZabcdefghijklmnopqrstuvwxyz
\font_fam_suf_fix
 \ssfam
 ss
 ABCDEFGHIJKLMNOPQRSTUVWXYZabcdefgijklmnopqrstuwxyz
\catcode`\_=8
\def\Cdss{{\fam\ssfam
    \mkern 4.2 mu \mathchoice%
    {\vrule height 6.5pt depth -.55pt width 1pt}%
    {\vrule height 6.5pt depth -.57pt width 1pt}%
    {\vrule height 4.55pt depth -.28pt width .8pt}%
    {\vrule height 3.25pt depth -.19pt width .6pt}%
    \mkern -6.3mu C}}%
\def\Fdss{{\fam\ssfam I\mkern -2.5mu F}}%
\def\Ndss{{\fam\ssfam I\mkern -2.5mu N}}%
\def\Qdss{{\fam\ssfam
    \mkern 3.8 mu \mathchoice%
    {\vrule height 6.5pt depth -.67pt width 1pt}%
    {\vrule height 6.5pt depth -.7pt width 1pt}%
    {\vrule height 4.55pt depth -.44pt width .7pt}%
    {\vrule height 3.25pt depth -.3pt width .5pt}%
    \mkern -5.9mu Q}}%
\def\Rdss{{\fam\ssfam I\mkern -2.5mu R}}%
\def\Zdss{{\fam\ssfam Z\mkern-8.1mu Z}}%
%
%
%
%
\font\teneuf=eufm10 
\font\seveneuf=eufm7
\font\fiveeuf=eufm5
\newfam\euffam \def\euf{\fam\euffam\teneuf} 
\textfont\euffam=\teneuf \scriptfont\euffam=\seveneuf
\scriptscriptfont\euffam=\fiveeuf

       \def\gfr{{\euf g}}

       \def\lfr{{\euf l}}

\input xy
\xyoption{matrix}
\xyoption{arrow}
\xyoption{curve}

\parindent=0pt

\def\Hom{{\rm Hom}}
\def\Ext{{\rm Ext}}
\def\bfb{\hbox{\bf b}}
\def\binom#1#2{\left({{#1}\atop{#2}}\right)}
\def\dlongrightarrow{\longrightarrow\hskip-8pt\rightarrow}

\centerline{\typd Algebras of p-adic distributions and}
\smallskip
\centerline{\typd admissible representations}

\medskip

\centerline{\typf P. Schneider, J. Teitelbaum}

\medskip

{\bf Introduction}

\smallskip

In a series of earlier papers, ([ST1-4]) we began a systematic
study of locally analytic representations of a locally
$L$-analytic group $G$, where $L\subseteq \Cdss_p$ is a finite
extension of $\Qdss_p$. Such a representation is given by a
continuous action of $G$ on a locally convex topological vector
space $V$ over a spherically complete extension field
$K\subseteq\Cdss_p$ of $L$, such that the orbit maps $g\mapsto gv$
are locally analytic functions on $G$. When $G$ is the group of
$L$-points of an algebraic group, the class of such
representations includes many interesting examples, such as the
principal series representations studied in [ST2], the finite
dimensional algebraic representations, and the smooth
representations of Langlands theory. A reasonable theory of such
representations requires the identification of a finiteness
condition that is broad enough to include the important examples
and yet restrictive enough to rule out pathologies. In this paper
we present such a finiteness condition that we call
"admissibility" for locally analytic representations (provided the
field $K$ is discretely valued). The admissible locally analytic
representations, among which are the examples mentioned above,
form an abelian category.

Our approach to the characterization of admissible representations
is based on the algebraic approach to such representations begun
in [ST2].  As in that paper, we require that the vector space $V$
carrying the locally analytic representation be of compact type, a
topological condition whose most important consequence is that $V$
is reflexive.   We focus our attention on  the algebra $D(G,K)$ of
locally analytic distributions on $G$. This algebra is the
continuous dual of the locally analytic, $K$-valued functions on
$G$, with multiplication given by convolution.  When $G$ is
compact, $D(G,K)$ is a Fr\'{e}chet algebra, but in general is
neither noetherian nor commutative.  If $V$ is a locally analytic
$G$-representation, then its continuous dual $V'_b$, with its
strong topology, becomes a module over $D(G,K)$. We identify a
subcategory of the module category of $D(G,K)$ that we call the
coadmissible modules. We show that any coadmissible  module
carries a canonical Fr\'echet topology.  We say that $V$ is
admissible if $V'_b$ is topologically isomorphic to a coadmissible
module. If $G$ is not compact, we say that $V$ is admissible if it
is admissible as a representation for one (or equivalently any)
compact open subgroup of $G$.

To define the category of coadmissible $D(G,K)$-modules, we
introduce a  general class of Fr\'{e}chet algebras that we call
Fr\'{e}chet-Stein algebras.  The definition of a Fr\'{e}chet-Stein
algebra, and its associated category of coadmissible modules,  is
inspired by the notion of a Stein algebra, the (commutative) ring
of functions on a Stein space, as described in  [For]. We recall
that for the commutative group $G = \Zdss_p^d$ the algebra
$D(G,\Qdss_p)$ is naturally isomorphic, via Fourier theory, to the
ring of rigid analytic functions $\Oscr(\Xscr)$ on the locally
analytic character group $\Xscr = \Xscr(\Zdss_p^d)$ of
$\Zdss_p^d$. This latter space $\Xscr$ is the $d$-dimensional open
unit disk over $\Qdss_p$. On a Stein space like $\Xscr$ the
category of coherent module sheaves is a well behaved abelian
category which, under the global section functor, is equivalent to
an abelian subcategory of the category ${\rm Mod}(\Oscr(\Xscr))$
of all $\Oscr(\Xscr)$-modules. A closer inspection of the proofs
of these facts reveals that, truly in the spirit of noncommutative
geometry, they are consequences of certain basic properties of the
ring $\Oscr(\Xscr)$ which can be formulated without making use of
the actual geometric object $\Xscr$ and which do not even require
the commutativity of the ring. From this observation we are led to
define a Fr\'echet-Stein algebra to be, roughly speaking, a
Fr\'echet algebra that is a projective limit of noetherian Banach
algebras, where the transition maps in the projective system are
flat.  A module over a Fr\'{e}chet-Stein algebra is coadmissible
if it is the projective limit of a compatible system of finitely
generated modules over the associated Banach algebras. The
category of coadmissible modules for a Fr\'{e}chet-Stein algebra
is abelian, and includes for example all finitely presented
$D(G,K)$-modules.

With the general notion of Fr\'{e}chet-Stein algebra in hand, our
first main result is that, when $G$ is compact, the algebra
$D(G,K)$ is Fr\'{e}chet-Stein. Our proof relies heavily on
Lazard's detailed study of $p$-adic analytic groups in his classic
work [Laz], and on the graded ring techniques presented in the
book [LVO].  The idea is that a sufficiently small $p$-valued
group is isomorphic as a $p$-adic manifold to $\Zdss_p^{d}$, and
so the vector space underlying the algebra $D(G,K)$ is isomorphic
to the dual of the $K$-valued locally analytic functions on
$\Zdss_p^{d}$. This space is known explicitly through the theory
of Mahler expansions. However, the multiplication on $D(G,K)$
comes from the non-commutative multiplication on $G$.  We use
graded techniques and Lazard's results to pass from commutative to
noncommutative information.

The idea to construct these Banach algebras appears for the first
time in the thesis of H. Frommer on the locally analytic principal
series of $GL_n$. He deals, by methods quite different from ours,
with the special case of the group $G$ of all unipotent triangular
matrices in $GL_n(\Zdss_p)$.

The Fr\'{e}chet-Stein property of $D(G,K)$ makes it possible to
define admissible representations as outlined above. Another
consequence is that any finitely
generated submodule of a finitely generated free $D(G,K)$-module
is closed.  It follows that any finitely presented $D(G,K)$-module
is analytic in the sense of [ST1], answering a question raised in
that paper.

Our second main result concerns the relation between the completed
group ring $\Zdss_p[[G]]$ of a compact $G$, studied in detail by
Lazard in [Laz], and the algebra $D(G,K)$.  In topological terms,
the ring $L[[G]]:=L\otimes \Zdss_p[[G]]$ is the Iwasawa algebra of
``measures'', and is  dual to the $L$-valued continuous functions
on $G$.  Because the locally analytic functions are a subspace of
the continuous functions, there is a natural map $L[[G]]\to
D(G,K)$. When $G$ is compact and locally $\Qdss_p$-analytic, we
prove that this map is faithfully flat, answering a question
raised in [ST3]. This result implies that, in an admissible
continuous representation of such a $G$ on a Banach space, as
studied in [ST3], the subspace of locally analytic vectors is
dense, and ``passage to analytic vectors'' is an exact functor.

The final topic of this paper is the development of a dimension
theory for coadmissible modules over $D(G,K)$. We begin by
considering such a theory for a general Fr\'{e}chet-Stein algebra
$A$. Our approach is modelled on the dimension theory of coherent
sheaves on Stein spaces as presented in [Ban]. In order to have a
well-defined "local" notion of dimension, we require that the
Banach algebras defining the Fr\'{e}chet-Stein structure of $A$ be
Auslander regular rings of bounded global dimension. Under this
assumption one can adopt the point of view that the "local
codimension of the support of a sheaf" is measured by the grade of
the module over the Banach algebra. We then show that
correspondingly the grade of a coadmissible module over $A$ has
the usual properties of a "global" codimension.

Having developed these general properties of codimension, we use
graded techniques and our earlier results to show that $D(G,K)$
satisfies the necessary regularity conditions when $G$ is a
compact, $d$-dimensional, $\Qdss_p$-analytic group. Consequently
every coadmissible module over $D(G,K)$ has a well-defined
codimension bounded above by $d$, as well as a dimension
filtration by coadmissible submodules. We also will show that, if
$V$ is a smooth admissible representation, or more generally a
$U(\gfr)$-finite admissible representation as in [ST1], then the
dual $D(G,K)$-module $V'_b$ has maximal codimension equal to $d$.

We begin the paper with a brief review of filtered and graded
techniques. After some preliminaries on Banach algebras, we define
Fr\'{e}chet-Stein algebras and their coadmissible modules and
establish the important properties of the category of such modules
in section 3.  Then we turn to $D(G,K)$, and after reviewing some
of Lazard's key results from [Laz] we prove the two main results
on $D(G,K)$ in sections 4 and 5. The application to representation
theory follows in section 6. Among other results we show that all
smooth admissible representations in the sense of Langlands theory
are admissible locally analytic representations. In fact the
former can be characterized among the latter by a trivial derived
action of the Lie algebra of $G$. In the final two sections we
discuss the existence of analytic vectors in continuous
representations, and develop the dimension theory of coadmissible
modules.

\medskip

{\bf Acknowledgements:} The second author was supported by a grant
from the US National Security Agency Mathematical Sciences
Program. Both authors gratefully acknowledge support from the
University of Illinois at Chicago and from the SFB ''Geometrische
Strukturen in der Mathematik'' at M\"unster making possible a
series of mutual visits.

\medskip

{\bf Notation:} For any associative unital ring $R$ we let
$\Mscr_R$ denote the full subcategory of all finitely generated
$R$-modules in the category ${\rm Mod}(R)$ of all (left) unital
$R$-modules. If $R$ is (left) noetherian then $\Mscr_R$ is an
abelian category.

Throughout $K$ is a field which is complete with respect to a
non-trivial nonarchimedean absolute value $|\ |$. For the basic
notions and facts in nonarchimedean functional analysis we refer to
[NFA].

Finally, $p$ is a fixed prime number.

\medskip

{\bf 1. Filtered rings}

\smallskip

If a ring is filtered then, quite often, it inherits properties of the
associated graded ring. Since this technique will be of fundamental
importance for us later on we begin by recalling some of it.

Throughout this section let $R$ be an associative unital ring. We call
$R$ {\it filtered} if it is equipped with a family
$(F^sR)_{s\in\Rdss}$ of additive subgroups $F^sR \subseteq R$ such
that, for any $r,s \in \Rdss$,

-- $F^rR \supseteq F^sR$ if $r \leq s$,

-- $F^rR \cdot F^sR \subseteq F^{r+s}R$,

-- $\bigcup_{s\in\Rdss} F^sR = R$ and $1 \in F^0R$.

For any $s \in \Rdss$ put
$$
F^{s+}R := \bigcup_{r > s} F^rR\ \ \ \ \hbox{and}\ \ \ \ gr^sR :=
F^sR/F^{s+}R\ .
$$
Then
$$
gr^\cdot R := \bigoplus_{s\in\Rdss} gr^sR
$$
with the obvious multiplication is called the {\it associated graded
ring}. The filtration is called {\it quasi-integral} if there exists
an $n_0 \in \Ndss$ such that $\{s\in\Rdss : gr^s R \neq 0\}
\subseteq \Zdss\cdot 1/n_0$. Usually we also assume that the filtered ring $R$ is
(separated and) {\it complete}, i.e., that the natural map
$$
R \mathop{\longrightarrow}\limits^{\cong}
\mathop{\lim\limits_{\longleftarrow}}\limits_{s} R/F^sR
$$
is bijective.

\medskip

{\bf Proposition 1.1:} {\it Let $R$ be a complete filtered ring whose
filtration is quasi-integral; if the graded ring $gr^\cdot R$ is
(left) noetherian then $R$ is (left) noetherian as well.}

Proof: Up to rescaling and reversing to increasing filtrations this is
[LVO] Prop. I.7.1.2.

\medskip

Any homomorphism $\phi : R \longrightarrow A$ between filtered rings
which respects the filtrations induces in the obvious way a
homomorphism of graded rings $gr^\cdot \phi : gr^\cdot R
\longrightarrow gr^\cdot A$.

In the following we consider two complete filtered rings $R$ and $A$
whose filtrations are quasi-integral and a unital ring homomorphism
$\phi : R \longrightarrow A$ which respects the filtrations.

\medskip

{\bf Proposition 1.2:} {\it  Suppose that $gr^\cdot R$ and $gr^\cdot
A$ are left noetherian and that $gr^\cdot A$ as a right $gr^\cdot
R$-module (via $gr^\cdot \phi$) is flat; then $A$ is flat as a right
$R$-module (via $\phi$).}

Proof: This is a combination of results in [LVO]. By Prop. II.1.2.1 it
suffices to check that

(a) good filtrations on filtered $A$-modules are separated and

(b) $F^\cdot R$ has the left Artin-Rees property.

According to Thm. I.5.7 the filtration $F^\cdot R$ induces a good
filtration on any left ideal in $R$. Hence, by Remark II.1.1.2(1), we
have property (b). Similarly, of course, $F^\cdot A$ has the
Artin-Rees property. Then, by Thm. II.1.1.5, the Rees ring
$\widetilde{A}$ is left noetherian which, by Lemma I.3.5.5(2) and Cor.
I.5.5(3), implies property (a).

\medskip

We will need a few intermediate steps in the above proofs as results
in their own right.

Let $J \subseteq R$ be a left ideal equipped with the filtration
induced by $F^\cdot R$. Similarly we equip the ideal $AJ$ in $A$ with
the filtration induced by $F^\cdot A$. Considering on $A \otimes_R J$
the tensor product filtration we then have the natural graded maps
$$
gr^\cdot A \otimes_{gr^\cdot R} gr^\cdot J \dlongrightarrow gr^\cdot
(A \otimes_R J) \longrightarrow gr^\cdot AJ
$$
of which the left one is surjective by construction.

\medskip

{\bf Lemma 1.3:} {\it  Suppose that $gr^\cdot R$ and $gr^\cdot A$ are
left noetherian and that $gr^\cdot A$ as a right $gr^\cdot R$-module
is flat; then }
$$
gr^\cdot A \otimes_{gr^\cdot R} gr^\cdot J
\mathop{\longrightarrow}\limits^{\cong} gr^\cdot (A
\otimes_R J) \mathop{\longrightarrow}\limits^{\cong} gr^\cdot AJ\ .
$$

Proof: By the flatness assumption the first and the third term in the
assertion both are submodules of $gr^\cdot A \otimes_{gr^\cdot R}\,
gr^\cdot R = gr^\cdot A$. It follows that $gr^\cdot A
\otimes_{gr^\cdot R}\, gr^\cdot J \cong gr^\cdot (A
\otimes_R J)$ and that $ gr^\cdot (A \otimes_R J) \hookrightarrow
gr^\cdot AJ$ is injective. Using [LVO] Cor. I.4.2.5 we may conclude
that the latter map even is an isomorphism provided the tensor product
filtration on $A \otimes_R J$ is separated. But this follows from
[LVO] Lemma I.6.6.15 and Prop. II.1.2.3 which, in particular, say that
the tensor product $A \otimes_R .$ preserves good filtrations and that
good filtrations (for $A$) are separated.

\medskip

{\bf Lemma 1.4:} {\it  Suppose that $gr^\cdot R$ is left noetherian;
there are nonzero elements $a_1,\ldots,a_l \in J$, say $a_i \in
F^{s_i}J\setminus F^{s_i+}J$, such that
$$
F^s J = \sum_{i=1}^l F^{s-s_i}R\cdot a_i
$$
for any $s \in \Rdss$. }

Proof: This is contained in the proof of [LVO] Thm. I.5.7.

\medskip

{\bf 2. Preliminaries on Banach algebras}

\smallskip

A $K$-Banach algebra $A$ is a $K$-Banach space $(A,|\ |_A)$ with a
structure of an associative unital $K$-algebra such that the
multiplication is continuous, which means that there is a constant $c
> 0$ such that
$$
|ab|_A \leq c |a|_A |b|_A\ \ \ \hbox{for any}\ a,b \in A\ .
$$
The norm will be called {\it submultiplicative} if it satisfies the
stronger condition that
$$
|1|_A = 1\ \ \ \hbox{and}\ \ \  |ab|_A \leq |a|_A |b|_A\ \ \
\hbox{for any}\ a,b \in A\ .
$$

\medskip

{\bf Proposition 2.1:} {\it Suppose that $A$ is a (left) noetherian
$K$-Banach algebra; we then have:

i. Each module $M$ in $\Mscr_A$ carries a unique $K$-Banach space
topology (called its $\underline{canonical}$ topology) such that the
$A$-module structure map $A \times M \rightarrow M$ is continuous;

ii. every $A$-submodule of a module in $\Mscr_A$ is closed in the
canonical topology; in particular, every (left) ideal in $A$ is
closed;

iii. any map in $\Mscr_A$ is continuous and strict for the canonical
topologies.}

Proof: The arguments in [BGR] 3.7.2 and 3.7.3 generalize in a
straightforward way to the noncommutative setting. We point out that,
by [BGR] Prop. 1.2.1/2, the given norm always can be replaced by an
equivalent one which is submultiplicative.

\medskip

In case the norm $|\ |_A$ is submultiplicative we may and will view
the $K$-Banach algebra $A$ as a complete filtered ring with respect to
the filtration
$$
F^s A := \{a \in A : |a|_A \leq p^{-s}\}\ .
$$
In general this filtration, of course, is not quasi-integral.

\medskip

{\bf 3. Fr\'echet-Stein algebras}

\smallskip

The following is very much inspired by the notions of a Stein
algebra and a Stein module in complex analysis as presented in
[For].

Let $A$ be a $K$-Fr\'echet algebra which means that $A$ is an
associative unital $K$-algebra such that the underlying $K$-vector
space has the structure of a $K$-Fr\'echet space and such that the
algebra multiplication is continuous. Consider a continuous seminorm
on $A$. It induces a norm on the quotient space $A/\{a\in A: q(a) =
0\}$. The completion of the latter with respect to $q$ is a $K$-Banach
space which will be denoted by $A_q$. It comes with a natural
continuous linear map $A \rightarrow A_q$ with dense image. For any
two continuous seminorms $q' \leq q$ the identity on $A$ extends to a
continuous, in fact norm decreasing, linear map $\phi^{q'}_q : A_q
\longrightarrow A_{q'}$ with dense image such that the diagram
$$
\xymatrix@R=0.5cm{
                &         A_q \ar[dd]^{\phi^{q'}_q}     \\
  A \ar[ur] \ar[dr]                 \\
                &         A_{q'}                 }
$$
commutes. For any sequence $q_1 \leq q_2 \leq\ldots\leq q_n
\leq\ldots$ of seminorms on $A$ which define the Fr\'echet topology
(such a sequence always exists), the obvious map
$$
A \mathop{\longrightarrow}\limits^{\cong}
\mathop{\lim\limits_{\longleftarrow}}\limits_{n\in \Ndss} A_{q_n}
$$
where on the right hand side the projective limit is formed with
respect to the $\phi^{q_n}_{q_{n+1}}$ as the transition maps is an
isomorphism of locally convex $K$-vector spaces. In the following a
continuous seminorm $q$ on $A$ will be called an algebra seminorm if
the multiplication on $A$ is continuous with respect to $q$, i.e., if
there is a constant $c > 0$ such that
$$
q(ab) \leq cq(a)q(b)\ \ \ \hbox{for any}\ a,b \in A\ .
$$
In this case $A_q$ in a natural way is a $K$-Banach algebra and the
natural map $A \rightarrow A_q$ is a homomorphism of $K$-algebras. If
the sequence $q_1 \leq\ldots\leq q_n \leq\ldots$ consists of algebra
seminorms then the transition maps $\phi^{q_n}_{q_{n+1}}$ are algebra
homomorphisms and
$$
A \mathop{\longrightarrow}\limits^{\cong}
\mathop{\lim\limits_{\longleftarrow}}\limits_{n\in \Ndss} A_{q_n}
$$
is an isomorphism of Fr\'echet algebras.

\medskip

{\bf Definition:} {\it The $K$-Fr\'echet algebra $A$ is called a
$K$-Fr\'echet-Stein algebra if there is a sequence $q_1 \leq\ldots\leq
q_n \leq\ldots$ of continuous algebra seminorms on $A$ which define
the Fr\'echet topology such that\hfill\break
 (i) $\ \ A_{q_n}$ is (left) noetherian, and\hfill\break
 (ii) $\ A_{q_n}$ is flat as a right $A_{q_{n+1}}$-module (via
 $\phi^{q_n}_{q_{n+1}}$)\hfill\break
for any $n\in\Ndss$.}


\medskip

Fix in the following a $K$-Fr\'echet-Stein algebra $A$ and a sequence
$(q_n)_{n\in\Ndss}$ as in the above definition. We want to introduce
the notions of coherent sheaves and coadmissible modules for $A$.

\medskip

{\bf Definition:} {\it A coherent sheaf for $(A,(q_n))$ is a family
$(M_n)_{n\in\Ndss}$ of modules $M_n$ in $\Mscr_{A_{q_n}}$ together
with isomorphisms $A_{q_n}
\otimes_{A_{q_{n+1}}} M_{n+1} \mathop{\longrightarrow}\limits^{\cong}
M_n$ in $\Mscr_{A_{q_n}}$ for any $n\in\Ndss$.}

\medskip

The coherent sheaves for $(A,(q_n))$ with the obvious notion of a
homomorphism form a category $Coh_{(A,(q_n))}$. As a consequence of
the flatness requirement in the definition of Fr\'echet-Stein this
category is abelian, again with the obvious notions of (co)kernels and
(co)images. For any coherent sheaf $(M_n)_n$ for $(A,(q_n))$ its
$A$-module of ``global sections'' is defined by
$$
\Gamma(M_n) := \mathop{\lim\limits_{\longleftarrow}}\limits_n
M_n\ .
$$

\medskip

{\bf Definition:} {\it A (left) $A$-module is called coadmissible if
it is isomorphic to the module of global sections of some coherent
sheaf for $(A,(q_n))$.}

\medskip

We let $\Cscr_A$ denote the full subcategory of coadmissible modules
in the category ${\rm Mod}(A)$. A simple cofinality argument shows
that $\Cscr_A$ indeed is independent of the choice of the sequence
$(q_n)_n$. Passing to global sections defines a functor
$$
\Gamma : Coh_{(A,(q_n))} \longrightarrow \Cscr_A\ .
$$

\medskip

{\bf Theorem:} {\it Let $(M_n)_n$ be a coherent sheaf for $(A,(q_n))$
and put $M := \Gamma(M_n)$; we have:

i. (Theorem A)  For any $n\in\Ndss$ the natural map $M
\longrightarrow M_n$ has dense image with respect to the canonical topology
on the target;

ii. (Theorem B)
$$
{\mathop{\lim\limits_{\longleftarrow}}\limits_n}^{(i)} M_n = 0
$$
for any natural number $i \geq 1$.}

Proof: i. Fix an $n\in\Ndss$. By [B-GT] II\S3.5 Thm. 1 it suffices to
show that, for any $m > n$, the image of $M_m$ in $M_n$ is dense.
Choose an $A_{q_m}$-linear surjection $A_{q_m}^r
\dlongrightarrow M_m$ from a finitely generated free
$A_{q_m}$-module onto $M_m$. Then all horizontal arrows in the
commutative diagram
$$
\xymatrix{
  A_{q_m}^r \ar[d] \ar@{->>}[r]
                & M_m \ar[d]  \\
  A_{q_n}^r \ar@{->>}[r]
                & A_{q_n}\otimes_{A_{q_m}} M_m = M_n }
$$
are surjective. Our claim now follows from the fact that the left
vertical arrow has dense image and that, by Prop. 2.1.iii, all the
maps in the diagram are continuous for the canonical topologies.

ii. The system $(M_n)_{n}$, by Theorem A, has the Mittag-Leffler
property as formulated in [EGA] III 0.13.2.4 so that loc. cit. 13.2.2
applies.

\medskip

It is immediate from Theorem B that the functor $\Gamma :
Coh_{(A,(q_n))}
\longrightarrow {\rm Mod}(A)$ is exact.

\medskip

{\bf Corollary 3.1:} {\it For any coherent sheaf $(M_n)_n$ for
$(A,(q_n))$ and $M := \Gamma(M_n)$ the natural map
$$
A_{q_n} \otimes_A M \mathop{\longrightarrow}\limits^{\cong} M_n
$$
is an isomorphism for any $n\in\Ndss$.}

Proof: By Theorem A the $A_{q_n}$-submodule of $M_n$ generated by the
image of $M$ is dense in $M_n$. Prop. 2.1.ii. then says that this
submodule, in fact, must be equal to $M_n$. This establishes the
surjectivity of the map in question and, more precisely, that $M_n$ as
an $A_{q_n}$ -module is generated by finitely many elements in the
image of the map $M \longrightarrow M_n$. Suppose now that $b_1\otimes
x_1 + \ldots + b_k\otimes x_k \in A_{q_n} \otimes_A M$ is an element
such that $b_1x_1 + \ldots + b_kx_k = 0$ in $M_n$. Consider the
homomorphism of coherent sheaves
$$
\matrix{
\hfill (A_{q_{n'}}^k)_{n'} & \longrightarrow & (M_{n'})_{n'} \hfill\cr\cr
(a_1,\ldots,a_k) & \longmapsto & a_1x_1 + \ldots a_kx_k\ . }
$$
Applying the above surjectivity argument to its kernel we find
finitely many elements $(c_1^{(1)},\ldots,c_k^{(1)}),\ldots,
(c_1^{(r)},\ldots,c_k^{(r)})$ in $A^k$ whose images generate the
kernel of the map $A_{q_n}^k \longrightarrow M_n$ as an
$A_{q_n}$-module. In particular there are $f_1,\ldots,f_r \in A_{q_n}$
such that
$$
(b_1,\ldots,b_k) = f_1\cdot (c_1^{(1)},\ldots,c_k^{(1)}) + \ldots +
f_r\cdot (c_1^{(r)},\ldots,c_k^{(r)})\ .
$$
It follows that
$$
\matrix{
\mathop{\sum}\limits_{i=1}^k b_i\otimes x_i & = &
\mathop{\sum}\limits_{i=1}^k \mathop{\sum}\limits_{j=1}^r f_j
c_i^{(j)}\otimes x_i = \mathop{\sum}\limits_{i=1}^k
\mathop{\sum}\limits_{j=1}^r f_j\otimes c_i^{(j)}x_i \cr\cr
& = & \mathop{\sum}\limits_{j=1}^r f_j\otimes
(\mathop{\sum}\limits_{i=1}^k c_i^{(j)}x_i) = 0\ .\hfill }
$$
Hence the map in question is injective as well.

\medskip

{\bf Remark 3.2:} {\it $A_{q_n}$, for any $n \in \Ndss$, is flat as a
right $A$-module.}

Proof: It suffices to show that for any left ideal $L \subseteq A$ the
induced map $A_{q_n} \otimes_A L \longrightarrow A_{q_n}$ is
injective. But this is exactly the same computation as in the previous
proof.

\medskip

{\bf Corollary 3.3:} {\it The functor $\Gamma : Coh_{(A,(q_n))}
\mathop{\longrightarrow}
\limits^{\sim} \Cscr_A$ is an equivalence of categories.}

Proof: By definition the functor is essentially surjective. According
to the previous corollary it is fully faithful. Both properties
together amount to the functor being an equivalence of categories.

\medskip

{\bf Corollary 3.4:} {\it i. The direct sum of two coadmissible
$A$-modules is coadmissible;

ii. the (co)kernel and (co)image of an arbitary $A$-linear map between
coadmissible $A$-modules are coadmissible;

iii. the sum of two coadmissible submodules of a coadmissible
$A$-module is coadmissible;

iv. any finitely generated submodule of a coadmissible $A$-module is
coadmissible;

v. any finitely presented $A$-module is coadmissible.}

Proof: The first assertion is obvious. The last three assertions are
immediate consequences of the first two. Hence it remains to establish
the second assertion. By the previous corollary any map between
coadmissible modules comes from a map between coherent sheaves. But,
by Theorem B, the functor $\Gamma$ into ${\rm Mod}(A)$ commutes with
the formation of (co)kernels and (co)images.

\medskip

{\bf Corollary 3.5:} {\it $\Cscr_A$ is an abelian subcategory of ${\rm
Mod}(A)$.}

\medskip

We want to discuss a few permanence properties of our notions. First
of all the subcategory $\Cscr_A$ in ${\rm Mod}(A)$ is not closed with
respect to the passage to submodules or quotient modules. To obtain a
precise statement we will equip any coadmissible $A$-module $M$ with a
topology. Write
$$
M = \mathop{\lim\limits_{\longleftarrow}}\limits_n M_n\ .
$$
By Prop. 2.1.i, each $M_n$ carries its canonical Banach space topology
as a finitely generated $A_{q_n}$-module. We equip $M$ with the
projective limit topology of these canonical topologies. This makes
$M$ into a $K$-Fr\'echet space. Moreover, the $A$-module structure map
$A \times M \rightarrow M$ clearly is continuous. We will call this
Fr\'echet topology the {\it canonical} topology of $M$.

\medskip

{\bf Lemma 3.6:} {\it For any coadmissible $A$-module $M$ and any
submodule $N \subseteq M$ the following assertions are equivalent:

i. $N$ is coadmissible;

ii. $M/N$ is coadmissible;

iii. $N$ is closed in the canonical topology of $M$.}

Proof: The equivalence of i. and ii. is immediate from Cor. 3.4.ii.
Suppose that $M = \Gamma(M_n)$. For any $n\in\Ndss$ let $N_n \subseteq
M_n$ denote the $A_{q_n}$-submodule generated by the image of $N$.
Then $(N_n)_n$ is a coherent subsheaf of $(M_n)_n$. Moreover, it is
straightforward to see that $\overline{N} :=
\Gamma(N_n)$ is the closure of $N$ in $M$ with its canonical
topology. Hence if $N$ is closed then it is coadmissible. On the other
hand, if $N$ is coadmissible then Cor. 3.1 implies that $N =
\Gamma(N_n) = \overline{N}$.

\medskip

It is immediate from Prop. 2.1.iii and Cor. 3.3 that any $A$-linear
map $f : M \longrightarrow N$ between two coadmissible $A$-modules is
continuous for the canonical topologies. It also is strict since its
image is a coadmissible submodule of $N$ by Cor. 3.4.ii and hence is
closed in $N$ by Lemma 3.6; so $f : M \longrightarrow {\rm im}(f)$ is
a continuous surjection between Fr\'echet spaces which by the open
mapping theorem is strict.

As another consequence of Cor. 3.4.iv and Lemma 3.6 we obtain that any
finitely generated left ideal of $A$ is closed.

\medskip

{\bf Proposition 3.7:} {\it Let $I$ be a closed two sided ideal in a
$K$-Fr\'echet-Stein algebra $A$; then $A/I$ is a $K$-Fr\'echet-Stein
algebra as well.}

Proof: By Lemma 3.6 both, $I$ and $A/I$, are coadmissible $A$-modules.
We therefore have the exact sequences
$$
0 \longrightarrow A_{q_n} \otimes_A I \longrightarrow A_{q_n}
\longrightarrow A_{q_n}/A_{q_n}I \longrightarrow 0
$$
and the isomorphism of $A$-modules
$$
A/I \mathop{\longrightarrow}\limits^{\cong}
\mathop{\lim\limits_{\longleftarrow}}\limits_n A_{q_n}/A_{q_n}I\ .
\leqno{(\ast)}
$$
Moreover $A_{q_n}I$ is the closure of the image of $I$ in $A_{q_n}$.
The latter implies that $A_{q_n}I$ is a two sided ideal in $A_{q_n}$.
Hence $A_{q_n}/A_{q_n}I$ is a $K$-Banach algebra with respect to the
quotient norm $\overline{q_n}$ of $q_n$. As a consequence of the open
mapping theorem $(\ast)$ then is an isomorphism of $K$-Fr\'echet
algebras. This means that $\overline{q_1} \leq\ldots\leq
\overline{q_n} \leq\ldots$ is a sequence of algebra seminorms on $A/I$ which
define the quotient topology. As a quotient of a noetherian algebra
$A_{q_n}/A_{q_n}I = (A/I)_{\overline{q_n}}$ is noetherian as well. By
construction we have
$$
\matrix{
(A/I)_{\overline{q_n}} & = A_{q_n} \otimes_A (A/I) = A_{q_n}
\otimes_{A_{q_{n+1}}} (A_{q_{n+1}} \otimes_A (A/I))\cr\cr
& = A_{q_n} \otimes_{A_{q_{n+1}}} (A/I)_{\overline{q_{n+1}}}\hfill }
$$
as $(A_{q_n},(A/I)_{\overline{q_{n+1}}})$-bimodules. The flatness of
$A_{q_n}$ as a right $A_{q_{n+1}}$-module therefore implies, by base
change (compare [B-CA] I\S2.7 Cor. 2), that $(A/I)_{\overline{q_n}}$
is flat as a right $(A/I)_{\overline{q_{n+1}}}$-module.

\medskip

{\bf Lemma 3.8:} {\it Let $A \longrightarrow B$ be a continuous unital
algebra homomorphism between $K$-Fr\'echet-Stein algebras such that
$B$ is coadmissible as a (left) $A$-module, and let $M$ be a (left)
$B$-module; then $M$ is coadmissible as a $B$-module if and only if it
is coadmissible as an $A$-module.}

Proof: Let $(q_n)_n$, resp. $(p_n)_n$, be a sequence of seminorms on
$A$, resp. $B$, as required in the definition of a Fr\'echet-Stein
algebra. Since the homomorphism $A \longrightarrow B$ is continuous we
find, for any $n\in\Ndss$, an $m \geq n$ such that this map extends
continuously to a homomorphism of Banach algebras $A_{q_m}
\longrightarrow B_{p_n}$. By passing to a subsequence of $(q_n)_n$
(and renumbering) we may assume that we always can choose $m = n$. We
then have, for any $n\in\Ndss$, the commutative diagram
$$
\xymatrix{
  A \ar[d] \ar[r] & B \ar[d]  \\
  A_{q_n} \ar[r] & B_{p_n}             }
$$
and hence the induced $(A_{q_n},B)$-bimodule map
$$
A_{q_n} \otimes_A B \longrightarrow B_{p_n}
$$
which is continuous for the canonical topology on the finitely
generated $A_{q_n}$-module $A_{q_n} \otimes_A B$. In the limit with
respect to $n$ these maps give back the identity map on $B$, but as a
continuous map from $B$ with its canonical topology as a coadmissible
$A$-module into $B$ with its given Fr\'echet topology. By the open
mapping theorem these two topologies then must in fact coincide.
Hence, for each $n\in\Ndss$, the canonical $(A,B)$-bimodule map $B
\longrightarrow A_{q_n} \otimes_A B$ is continuous and factorizes
therefore through a continuous $(A_{q_{m(n)}},B)$-bimodule map
$B_{p_{m(n)}} \longrightarrow A_{q_n} \otimes_A B$ for some $m(n) \geq
n$. This proves that we have
$$
\mathop{\lim\limits_{\longleftarrow}}\limits_n B_{p_n} \otimes_B M =
\mathop{\lim\limits_{\longleftarrow}}\limits_n (A_{q_n} \otimes_A B)
\otimes_B M =
\mathop{\lim\limits_{\longleftarrow}}\limits_n A_{q_n} \otimes_A M\ .
$$
If $A_{q_n} \otimes_A M$ is a finitely generated $A_{q_n}$-module then
$B_{p_n} \otimes_{A_{q_n}} (A_{q_n} \otimes_A M) = B_{p_n} \otimes_A
M$ and hence its quotient $B_{p_n} \otimes_B M$ are finitely generated
$B_{p_n}$-modules. Let us assume, vice versa, that $B_{p_{m(n)}}
\otimes_B M$ is a finitely generated $B_{p_{m(n)}}$-module. We choose
a set of generators $1 \otimes x_1,\ldots,1 \otimes x_k$ with $x_i
\in M$. We also fix generators $1 \otimes y_1,\ldots,1 \otimes y_l$,
with $y_j \in B$, of the finitely generated $A_{q_n}$-module $A_{q_n}
\otimes_A B$. Given any $x \in M$ we find $b_i \in B_{p_{m(n)}}$ such
that
$$
1 \otimes x = \sum_{i=1}^k b_i \otimes x_i \leqno{(1)}
$$
in $B_{p_{m(n)}} \otimes_B M$. For each $b_i$ we then find $c_{i,j}
\in A_{q_n}$ such that
$$
\sum_{j=1}^l c_{i,j} \otimes y_j = \hbox{image of}\ b_i\ \hbox{in}\
A_{q_n} \otimes_A B\ .\leqno{(2)}
$$
By first applying the map $B_{p_{m(n)}} \otimes_B M \longrightarrow
(A_{q_n} \otimes_A B) \otimes_B M = A_{q_n} \otimes_A M$ and then
inserting $(2)$ the identity $(1)$ becomes
$$
1 \otimes 1 \otimes x = \sum_{i,j} c_{i,j} \otimes y_j \otimes x_i,\ \
\hbox{resp.}\ \ 1 \otimes x = \sum_{i,j} c_{i,j} \otimes y_jx_i\ .
$$
This proves that the finitely many elements $1 \otimes y_jx_i$
generate $A_{q_n} \otimes_A M$ as an $A_{q_n}$-module.

\medskip

In the setting of Prop. 3.7 this last lemma means that
$$
\Cscr_{A/I} = {\rm Mod}(A/I) \cap \Cscr_A\ .
$$
We finish this section with a general simplicity criterion for
coadmissible modules.

\medskip

{\bf Lemma 3.9:} {\it Let $M$ be a coadmissible $A$-module
corresponding to the coherent sheaf $(M_n)_n$; the $A$-module $M$ is
simple if the $A_{q_n}$-module $M_n$ is simple for infinitely many
$n\in\Ndss$.}

Proof: Suppose $M$ is not simple and is nonzero. We then find a
nonzero finitely generated proper $A$-submodule $0 \neq N \subset M$.
By Cor. 3.4.iv $N$ is coadmissible and therefore corresponds to a
nonzero proper subsheaf $0 \neq (N_n)_n \subset (M_n)_n$. There have
to exist $n_0, n_1 \in \Ndss$ such that $N_{n_0} \neq 0$ and $N_{n_1}
\neq M_{n_1}$. As a consequence of Cor. 3.1 any $n \geq
\sup(n_0,n_1)$ then satisfies $N_n \neq 0$ and $N_n \neq M_n$. In
other words, for these $n$ the $A_{q_n}$-module $M_n$ is not simple.


\medskip

{\bf 4. Distribution algebras of uniform pro-$p$-groups}

\smallskip

Let $\Qdss_p \subseteq L \subseteq K \subseteq
\Cdss_p$ be complete intermediate fields where $L/\Qdss_p$ is
finite and $K$ is discretely valued. The absolute value $|\ |$ on
$\Cdss_p$ is normalized as usual by $|p| = p^{-1}$. Furthermore,
let $G$ be a compact locally $L$-analytic group. In [ST2] \S2 we
have introduced the $K$-Fr\'echet algebra $D(G,K)$ of $K$-valued
locally analytic distributions on $G$. One of the principal aims
of this paper is to show that $D(G,K)$ has a natural structure of
a $K$-Fr\'echet-Stein algebra. The crucial case to consider is $L
= \Qdss_p$ and $G$ a locally $\Qdss_p$-analytic group which in
addition is a uniform pro-$p$-group. We remark right away that any
compact locally $\Qdss_p$-analytic group has an open normal
subgroup which is a uniform pro-$p$-group (compare [DDMS] Cor.
8.34). The technical heart of the matter is the use and
appropriate generalization of the methods of Lazard in [Laz]. It
therefore seems quite natural and, in fact, allows for more
flexibility to also use his language. So we begin by reviewing
some of his basic notions and results.

A {\it $p$-valuation} on a group $G$ is a real valued function $\omega
: G\setminus\{1\} \longrightarrow (1/(p-1),\infty)$ such that
$$
\matrix{
\omega(gh^{-1}) \geq {\rm min}(\omega(g),\omega(h)),\hfill\cr
\omega(g^{-1}h^{-1}gh) \geq \omega(g) + \omega(h), \hbox{and}\hfill\cr
\omega(g^p)=\omega(g) + 1\hfill }
$$
for any $g,h \in G$ ([Laz] III.2.1.2). As usual one puts $\omega(1) :=
\infty$. For each real number $\nu > 0$ we define the normal
subgroups
$$
G_\nu := \{ g \in G : \omega(g) \geq \nu\}\ \ \ \hbox{and}\ \ \
G_{\nu^+} := \{ g \in G : \omega(g) > \nu\}
$$
of $G$, and we put
$$
gr(G) := \bigoplus_{\nu > 0} G_\nu/G_{\nu^+}\ .
$$
With the Lie bracket induced by the commutator $gr(G)$ is a graded Lie
algebra over $\Fdss_p$. Let $\Fdss_p[\epsilon]$ denote the polynomial
ring in one variable $\epsilon$ over $\Fdss_p$, as usual viewed as a
graded $\Fdss_p$-algebra with $\epsilon$ of degree $1$. The rule
$gG_{\nu^+} \mapsto g^pG_{(\nu +1)^+}$ defines an $\Fdss_p$-linear
operator $P$ on $gr(G)$, which is homogeneous of degree $1$, and which
satisfies $[P\overline{g},\overline{h}] =
P([\overline{g},\overline{h}])$ for homogeneous elements
$\overline{g},\overline{h}$ of $gr(G)$. Letting $\epsilon$ act as $P$
therefore makes $gr(G)$ into a graded Lie algebra over
$\Fdss_p[\epsilon]$ ([Laz] III.2.1.1). As an
$\Fdss_p[\epsilon]$-module $gr(G)$ is free and its rank is called the
{\it rank} of $G$ ([Laz] III.2.1.3).

For the remainder of this section we fix a compact locally
$\Qdss_p$-analytic group $G$ together with a $p$-valuation
$\omega$ on $G$. As a consequence of [Laz] III. 3.1.3/7/9 we have:

-- The topology of $G$ is defined by $\omega$ ([Laz] II.1.1.5); in
particular $G$ is a pro-$p$-group.

-- The rank of $G$ is equal to the dimension $d$ of $G$ (as a locally
analytic manifold) and in particular is finite.

By [Laz] III.3.1.11 a given $p$-valuation on $G$ can always be changed
into one which has rational values. We therefore assume from now on
that $\omega$ has rational values.

An {\it ordered basis} of $G$ is a sequence of elements
$h_1,\ldots,h_d\in G\setminus\{1\}$ such that the elements
$h_iG_{\omega(h_i)^+} \in gr(G)$ form a basis of $gr(G)$ as an
$\Fdss_p[\epsilon]$-module. It always exists. Given one the map
$$
\matrix{
\psi : &\Zdss_p^d & \mathop{\longrightarrow}\limits^\sim & G \cr\cr
& (x_1,\ldots x_d) & \longmapsto & h_1^{x_1}\cdot\ldots\cdot
h_d^{x_d}}
$$
is a bijective global chart for the manifold $G$ satisfying
$$
\omega(h_1^{x_1}\cdot\ldots\cdot h_d^{x_d}) = \mathop{\rm
min}\limits_{1 \leq i \leq d} (\omega(h_i) + \omega_p(x_i))
$$
where $\omega_p$ denotes the $p$-adic valuation on $\Zdss_p$
([Laz] III.2.2.5/6). In the following we fix an ordered basis
$(h_1,\ldots,h_d)$ of $G$ together with the corresponding chart
$\psi$. Using this chart we may identify the locally convex
$K$-vector spaces of locally analytic functions
$$
\psi^{\ast} : C^{an}(G,K) \mathop{\longrightarrow}\limits^{\cong}
C^{an}(\Zdss_p^d,K)
$$
as well as the larger $K$-Banach spaces of continuous functions
$$
\psi^{\ast} : C(G,K) \mathop{\longrightarrow}\limits^{\cong}
C(\Zdss_p^d,K)
$$
on both sides. Via Mahler expansions (compare [Laz] III.1.2.4)
$C(\Zdss_p^d,K)$ can be viewed as the space of all series
$$
f(x) =
\sum_{\alpha\in\Ndss_0^{d}} c_{\alpha}\binom{x}{\alpha}
$$
with $c_\alpha \in K$ and such that $|c_\alpha| \longrightarrow 0$ as
$|\alpha| \longrightarrow \infty$. Here we put, as usual,
$$
\binom{x}{\alpha} := \binom{x_1}{\alpha_1}\cdots\binom{x_d}{\alpha_d}
$$
and
$$
|\alpha| := \sum_{i=1}^{d} \alpha_{i}
$$
for $x = (x_1,\ldots,x_d) \in \Zdss_p^d$ and multi-indices $\alpha =
(\alpha_1,\ldots,\alpha_d) \in \Ndss_0^d$. By Amice's theorem (compare
[Laz] III.1.3.9) the Mahler expansion $f$ lies in the subspace
$C^{an}(\Zdss_p^d,K)$ if and only if $|c_\alpha|r^{|\alpha|}
\longrightarrow 0$ for some real number $r > 1$ as $|\alpha|
\longrightarrow \infty$.

We embed the group ring $\Zdss_p[G]$ into the distribution algebra
$D(G,K)$ by viewing a group element $g \in G$ as the Dirac
distribution $\delta_g$. We then have
$$
g(f) = \delta_{\psi(x)}(f) = \psi^{\ast}(f)(x)
$$
for any $f \in C(G,K)$ and any $g = \psi(x) \in G$. Henceforth we
write $b_i := h_i - 1$ and $\bfb^{\alpha} :=
b_1^{\alpha_1}b_2^{\alpha_2}\cdots b_{d}^{\alpha_d}$, for $\alpha
= (\alpha_1,\ldots,\alpha_d) \in \Ndss_0^d$, in $\Zdss_p[G]
\subseteq D(G,K)$. If $c_\alpha$ denote the coefficients of the
Mahler expansion of $\psi^{\ast}(f)$ for some $f \in C(G,K)$ then
their classical computation in terms of finite differences shows
that
$$
\bfb^{\alpha}(f)=c_\alpha\ .
$$
It easily follows that any distribution $\lambda \in D(G,K)$ has a
unique convergent expansion
$$
\lambda = \sum_{\alpha\in\Ndss_0^{d}} d_{\alpha}\bfb^{\alpha}
$$
with $d_\alpha \in K$ such that, for any $0 < r < 1$, the set
$\{|d_\alpha|r^{|\alpha|}\}_{\alpha\in\Ndss_0^d}$ is bounded.
Conversely, any such series is convergent in $D(G,K)$. The Fr\'echet
topology on $D(G,K)$ is defined by the family of norms
$$
\|\lambda\|'_r := \mathop{\rm sup}\limits_{\alpha\in\Ndss_0^d}
|d_\alpha|r^{|\alpha|}
$$
for $0 < r < 1$. For our purposes it is in fact more convenient to
take the original $p$-valuation $\omega$ into account and to always
work with the equivalent set of norms
$$
\|\lambda\|_r := \mathop{\rm sup}\limits_{\alpha\in\Ndss_0^d}
|d_\alpha|r^{\tau\alpha}
$$
where $\tau\alpha := \sum_i \alpha_i\omega(h_i)$. We note that, since
the multiplication in $D(G,K)$ is jointly continuous, we obtain the
expansion of the product of two distributions by multiplying their
expansions, inserting the expansions
$$
\bfb^\beta \bfb^\gamma = \sum_\alpha
c_{\beta\gamma,\alpha}\bfb^\alpha\ ,
$$
and rearranging.

The inclusion $\Zdss_p[G] \subseteq D(G,K)$ extends to an embedding of
topological rings
$$
\Zdss_p[[G]] \hookrightarrow D(G,K)
$$
of the completed group ring
$$
\Zdss_p[[G]] := \mathop{\rm lim}\limits_{\longleftarrow} \Zdss_p[G/N]
$$
with $N$ running over all open normal subgroups of $G$ into the
distribution algebra. A distribution $\lambda =
\sum_{\alpha} d_{\alpha}\bfb^{\alpha}$ lies in $\Zdss_p[[G]]$ if and
only if all $d_\alpha \in \Zdss_p$; moreover, the norm $\|\ \|_{1/p}$
restricted to $\Zdss_p[[G]]$ is independent of the choice of the
ordered basis of $G$, is multiplicative, satisfies $\|g-1\|_{1/p} \leq
p^{-\omega(g)}$ for any $g \in G$, and defines the compact topology of
$\Zdss_p[[G]]$ ([Laz] III.2.3; his $w$ is the additive version of our
$\|\ \|_{1/p}$). The latter in fact implies that each of the norms
$\|\ \|_r$, for $0 < r < 1$, defines this same compact topology of
$\Zdss_p[[G]]$. Theorem III.2.3.3 in [Laz] says that the inclusion map
$G
\hookrightarrow \Zdss_p[[G]]$ induces an isomorphism
$$
U(gr(G)) \mathop{\longrightarrow}\limits^{\cong} gr^\cdot_{1/p}
\Zdss_p[[G]]
$$
where the left hand side is the universal enveloping algebra of the
Lie algebra $gr(G)$ over $\Fdss_p[\epsilon]$ and the right hand side
is the associated graded ring for the filtration
$$
F^s_{1/p} \Zdss_p[[G]] := \{ \lambda \in \Zdss_p[[G]] :
\|\lambda\|_{1/p} \leq p^{-s}\}\ .
$$
Note that the right hand side is naturally an
$\Fdss_p[\epsilon]$-algebra as well by viewing $\Fdss_p[\epsilon]$ as
the associated graded ring for the $p$-adic valuation on $\Zdss_p$.

After this review we begin our investigation by looking at the
expansions
$$
\bfb^\beta \bfb^\gamma = \sum_\alpha
c_{\beta\gamma,\alpha}\bfb^\alpha
$$
with $c_{\beta\gamma,\alpha} \in \Zdss_p$.

\medskip

{\bf Lemma 4.1:} {\it i. $\omega_p(c_{\beta\gamma,\alpha}) \geq {\rm
max}(0,\tau\beta + \tau\gamma - \tau\alpha)$;

ii. $|c_{\beta\gamma,\alpha}|r^{\tau\alpha} \leq  r^{\tau\beta +
\tau\gamma}$ for any $1/p \leq r \leq 1$.}

Proof: Since $\|\ \|_{1/p}$ is multiplicative on $\Zdss_p[[G]]$ we
have
$$
p^{-\tau\beta - \tau\gamma} = \|\bfb^{\beta}\bfb^{\gamma}\|_{1/p} =
\mathop{\rm sup}\limits_{\alpha} |c_{\beta\gamma,\alpha}|p^{-\tau\alpha}
$$
which implies $\omega_p(c_{\beta\gamma,\alpha}) \geq \tau\beta +
\tau\gamma - \tau\alpha$. It also implies the second assertion in case
$r = 1/p$. For general $r$ we distinguish two cases. If $\tau\alpha
\geq \tau\beta + \tau\gamma$ then $|c_{\beta\gamma,\alpha}| \leq 1
\leq r^{\tau\beta + \tau\gamma - \tau\alpha}$. If $\tau\alpha
\leq \tau\beta + \tau\gamma$ then $|c_{\beta\gamma,\alpha}| \leq
(1/p)^{\tau\beta + \tau\gamma - \tau\alpha} \leq r^{\tau\beta +
\tau\gamma - \tau\alpha}$.

\medskip

{\bf Proposition 4.2:} {\it Each norm $\|\ \|_r$ on $D(G,K)$, for $1/p
\leq r < 1$, is submultiplicative.}

Proof: Let $\lambda = \sum_{\beta} d_{\beta}\bfb^{\beta}$ and $\mu =
\sum_{\gamma} e_{\gamma}\bfb^{\gamma}$ be two distributions in
$D(G,K)$. Their product then has the expansion
$$
\lambda\mu = \sum_{\beta,\gamma}
d_{\beta}e_{\gamma}\bfb^{\beta}\bfb^{\gamma} = \sum_{\alpha}
(\sum_{\beta,\gamma} d_{\beta}e_{\gamma}c_{\beta\gamma,\alpha})
\bfb^{\alpha}\ .
$$
Using Lemma 4.1.ii we obtain
$$
\matrix{
\|\lambda\mu\|_r & = \mathop{\rm sup}\limits_{\alpha}
|\mathop{\sum}\limits_{\beta,\gamma}
d_{\beta}e_{\gamma}c_{\beta\gamma,\alpha}| r^{\tau\alpha} \leq
\mathop{\rm sup}\limits_{\alpha,\beta,\gamma}
|d_{\beta}||e_{\gamma}||c_{\beta\gamma,\alpha}|r^{\tau\alpha}
\hfill\cr\cr
& \leq \mathop{\rm sup}\limits_{\beta,\gamma}
|d_{\beta}|r^{\tau\beta}|e_{\gamma}|r^{\tau\gamma} \hfill\cr\cr &
\leq \|\lambda\|_r \|\mu\|_r\ .\hfill }
$$

\medskip

We put
$$
K[[G]] := K \otimes_{\Qdss_p} \Zdss_p[[G]]
$$
which still is a subalgebra of $D(G,K)$. On any of the rings $R =
\Zdss_p[[G]], K[[G]]$, or $D(G,K)$ we have, for each $1/p \leq r < 1$, the
filtration
$$
F^s_r R := \{a \in R : \|a\|_r \leq p^{-s}\}\ .
$$
The associated graded ring is denoted by $gr^{\cdot}_r R$. We let
$$
D_r(G,K) := \hbox{completion of}\ D(G,K)\ \hbox{with respect to the
norm}\ \|\ \|_r\ .
$$
As a $K$-Banach space $D_r(G,K)$ is given by all series
$$
\sum_\alpha d_\alpha \bfb^\alpha\ \ \hbox{with}\ d_\alpha \in K\
\hbox{and}\ |d_\alpha|r^{\tau\alpha} \longrightarrow 0\ \hbox{as}\
|\alpha| \longrightarrow \infty\ .
$$
We also introduce the even larger $K$-Banach space
$$
D_{<r}(G,K) := \{\sum_\alpha d_\alpha \bfb^\alpha : d_\alpha \in K\
\hbox{and}\ \{|d_\alpha|r^{\tau\alpha}\}_\alpha\ \hbox{is bounded}\}\
.
$$
On both of them the norm continues to be given by
$$
\|\sum_\alpha d_\alpha \bfb^\alpha\|_r =  \mathop{\rm sup}\limits_{\alpha}
|d_\alpha|r^{\tau\alpha}\ .
$$
By Prop. 4.2 the multiplication on $D(G,K)$ extends continuously to
$D_r(G,K)$ and makes $D_r(G,K)$ into a $K$-Banach algebra. If $r >
1/p$ then $D_{<r}(G,K) \subseteq D_{1/p}(G,K)$. The computation in the
proof of Prop. 4.2 shows that the left term is multiplicatively closed
in the algebra $D_{1/p}(G,K)$ and that $D_{<r}(G,K)$ with the norm
$\|\ \|_r$ is a $K$-Banach algebra. Altogether we obtain the system of
$K$-Banach algebras
$$
\ldots\subseteq D_{r}(G,K) \subseteq D_{<r}(G,K) \subseteq D_{r'}(G,K)
\subseteq\ldots\subseteq D_{1/p}(G,K)
$$
with $1/p \leq r' < r < 1$ and
$$
D(G,K) = \mathop{\lim\limits_{\longleftarrow}}\limits_r D_r(G,K) =
\mathop{\lim\limits_{\longleftarrow}}\limits_r D_{<r}(G,K)\ .
$$
On $R = D_r(G,K)$, resp. $= D_{<r}(G,K)$, we again have, for any $1/p
\leq r' \leq r$, resp. $1/p \leq r' < r$, the filtration
$$
F^s_{r'} R := \{a \in R : \|a\|_{r'} \leq p^{-s}\}
$$
with associated graded ring $gr^\cdot_{r'} R$.

It is crucial for our results that all the various graded rings we
have introduced can be computed explicitly. First of all we remark
that all the norms under consideration restrict on $\Zdss_p$ to the
$p$-adic absolute value. Hence all graded rings naturally are
$\Fdss_p[\epsilon]$-algebras. Secondly we make the observation that,
since $K[[G]]$ is dense in $D_r(G,K)$, we have the isomorphism
$$
gr^\cdot_r K[[G]] \mathop{\longrightarrow}\limits^{\cong} gr^\cdot_r
D_r(G,K)
$$
for any $1/p \leq r < 1$. Moreover, by construction, the natural map
of $\Fdss_p[\epsilon]$-algebras $ gr^\cdot_r \Zdss_p[[G]]
\hookrightarrow gr^\cdot_r K[[G]]$ is injective. It extends to a ring
homomorphism
$$
gr^\cdot K \mathop{\otimes}\limits_{\Fdss_p[\epsilon]} gr^\cdot_r
\Zdss_p[[G]]
\longrightarrow gr^\cdot_r K[[G]]
$$
where $gr^\cdot K$ is the associated graded ring for the
filtration $F^sK := \{ a \in K : |a| \leq p^{-s}\}$. This map, in
fact, is an isomorphism since $gr^\cdot_r \Zdss_p[[G]]$, resp.
$gr^\cdot_r K[[G]]$, is free as a module over $\Fdss_p[\epsilon]$,
resp. $gr^\cdot K$, with basis the principal symbols
$\sigma_r(\bfb^\alpha)$. We recall that the {\it principal symbol}
$\sigma_r(a)$ of a nonzero element $a$ in one of our rings $R$
with respect to the filtration $F_r^\cdot R$ is the coset
$\sigma_r(a) := a + F_r^{s+}R$ if $a \in F_r^s R \setminus
F_r^{s+}R$.

\medskip

{\bf Lemma 4.3:} {\it For any $1/p \leq r < 1$ we have}
$$
gr^\cdot K \mathop{\otimes}\limits_{\Fdss_p[\epsilon]} gr^\cdot_r
\Zdss_p[[G]]
\mathop{\longrightarrow}\limits^{\cong} gr^\cdot_r K[[G]]
\mathop{\longrightarrow}\limits^{\cong} gr^\cdot_r D_r(G,K)
\ .
$$

\medskip

We note that
$$
gr^\cdot K \cong k[\epsilon_{\rm o},\epsilon_{\rm o}^{-1}]
$$
where $k$ is the residue class field of $K$ and $\epsilon_{\rm o}$
denotes the principal symbol of a prime element of $K$. The natural
map $\Fdss_p[\epsilon] \hookrightarrow gr^\cdot K$ is given by
$\epsilon
\mapsto  a\epsilon_{\rm o}^e$ where $e$ denotes the absolute ramification index of $K$
and $a \in \Fdss_p^\times$ is an appropriate element.

To obtain our results we have to impose additional condition on our
$p$-valued group $(G,\omega)$. According to [Laz] III.2.1.6 the pair
$(G,\omega)$ is called {\it $p$-saturated} if any $g \in G$ such that
$\omega(g) > p/(p-1)$ is a $p^{\rm th}$ power. In the following we
make the assumption that
$$
\matrix{
(G,\omega)\ \hbox{is $p$-saturated and the ordered basis}\
(h_1,\ldots,h_d)\ \hbox{of}\ G\cr  \hbox{satisfies}\ \omega(h_i) +
\omega(h_j) > p/(p-1)\ \hbox{for any}\ 1 \leq i \neq j \leq d\ .\hfill }
\leqno{({\rm HYP})}
$$

\medskip

{\bf Remark:} In [DDMS] \S4.1 the authors introduce the class of
uniform (or uniformly powerful) pro-$p$-groups. Without recalling
their definition we briefly relate this notion to the language of
Lazard which we are using. For simplicity we assume that $p \neq
2$. First let $G$ be, as above, a compact locally
$\Qdss_p$-analytic group which has a $p$-valuation $\omega$ such
that (HYP) holds. As we have seen in the proof of Lemma 4.4 the
commutator of any two members of the given ordered basis is a
$p$-th power. This easily implies that $G$ is powerful ([DDMS]
Def. 3.1). On the other hand any group with a $p$-valuation must
be torsionfree. According to [DDMS] Thm. 4.5 these two properties
together ensure that $G$ is uniform.\hfill\break
 Vice versa, assume $G$ to be a uniform pro-$p$-group. It is locally
$\Qdss_p$-analytic by [DDMS] Thm. 8.18. The discussion in [DDMS]
\S4.2 in fact shows that $G$ has an integrally valued
$p$-valuation $\omega'$ and an ordered basis $(h_1,\ldots,h_d)$
such that $\omega'(h_1) = \ldots = \omega'(h_d) = 1$; moreover, by
[DDMS] Lemma 4.10, each element in $G_2$ is a $p$-th power.
(Compare also the Notes at the end of Chap. 4 in [DDMS].) Hence
(HYP) holds in this context.

\medskip

This remark in particular shows that we could have worked from the
beginning with such simpler $p$-valuations like $\omega'$. But in
order to have as much flexibility as possible for later
applications it seemed advantageous  to us to instead minimize the
requirements on the $p$-valuation.

\medskip

{\bf Lemma 4.4:} {\it Assuming (HYP) we have $\|b_ib_j - b_jb_i\|_r <
\|b_ib_j\|_r$ for any $1 \leq i < j \leq d$ and any $1/p < r < 1$.}

Proof: Put $c := \omega(h_i) + \omega(h_j)$. Since $i < j$ we have
$\|b_ib_j\|_r = r^c$. To estimate the left hand side in the assertion
we use that our hypothesis implies that $h := h_i^{-1}h_j^{-1}h_ih_j =
g^p$ for some $g \in G$. Hence
$$
\matrix{
b_ib_j - b_jb_i & = h_ih_j - h_jh_i = h_jh_i(h - 1) \hfill\cr\cr
                & = h_jh_i(g^p - 1) = h_jh_i(((g-1) + 1)^p - 1) \cr\cr
                & = h_jh_i(g-1)^p + \mathop{\sum}\limits_{n=1}^{p-1}
                    \pmatrix{p \cr n} h_jh_i(g-1)^n \hfill }
$$
and therefore, by the submultiplicativity of $\|\ \|_r$ and the fact
that $\|h_i\|_r = 1$,
$$
\|b_ib_j - b_jb_i\|_r \leq {\rm max}(\|g-1\|_r^p,|p|\|g-1\|_r)\ .
$$
Since $1/p \leq r < 1$ the inequality $\|g-1\|_{1/p} \leq
p^{-\omega(g)}$ implies that $\|g-1\|_r \leq r^{\omega(g)}$. To see
this let $0 < z \leq 1$ such that $(1/p)^z = r$ and consider the
expansion $g-1 = \sum_\alpha d_\alpha\bfb^\alpha$. We have
$|d_\alpha|p^{-\tau\alpha} \leq p^{-\omega(g)}$. By exponentiating
this latter inequality and using that $|d_\alpha| \leq 1$ we have
$$
|d_\alpha|r^{\tau\alpha} \leq |d_\alpha|^z p^{-z\tau\alpha} \leq
p^{-z\omega(g)} = r^{\omega(g)}\ .
$$
Combining this with the previous estimate we obtain
$$
\|b_ib_j - b_jb_i\|_r \leq {\rm
max}(r^{p\omega(g)},p^{-1}r^{\omega(g)})\ .
$$
>From the properties of a $p$-valuation we deduce that
$$
\omega(g) = \omega(h) - 1 \geq \omega(h_i) + \omega(h_j) - 1 = c - 1\
.
$$
We therefore get
$$
\|b_ib_j - b_jb_i\|_r \leq {\rm
max}(r^{p(c-1)},p^{-1}r^{c-1})\ .
$$
The inequality $p^{-1}r^{c-1} < r^{c}$ is obvious from $1/p < r$. It
remains to be seen that $c < p(c-1)$. But this is equivalent to
$p/(p-1) < c$ which is part of our assumptions.

\medskip

{\bf Theorem 4.5:} {\it Assuming (HYP) let $1/p < r < 1$ and let $R$
be one of the rings $\Zdss_p[[G]]$ or $D_r(G,K)$; we then have:

i. $gr^\cdot_r R$ is a polynomial ring over $\Fdss_p[\epsilon]$, resp.
$gr^\cdot K$, in the variables
$\sigma_r(b_1),\ldots,\break\sigma_r(b_d)$; the norm $\|\ \|_r$ on $R$
is multiplicative;

ii. if $r \in p^\Qdss$ then $R$ is a (left and right) noetherian
integral domain. }

Proof: i. By Lemma 4.3 we only need to consider the case $R =
\Zdss_p[[G]]$. We have remarked already that $gr^\cdot_r \Zdss_p[[G]]$
as a module is free over $\Fdss_p[\epsilon]$ with basis
$\{\sigma_r(\bfb^\alpha )\}_\alpha$. Since $\|\bfb^\alpha\|_r =
\prod_i \|b_i\|_r^{\alpha_i}$ we have $\sigma_r(\bfb^\alpha) = \prod_i
\sigma_r(b_i)^{\alpha_i}$. Moreover, according to Lemma 4.4, the
$\sigma_r(b_i)$ commute with each other. The norm $\|\ \|_r$ has to be
multiplicative since $gr^\cdot_r R$ is an integral domain.

ii. Because of our assumptions that $K$ is discretely valued, $\omega$
has rational values, and $r \in p^\Qdss$ the filtration $F^\cdot_r R$
is quasi-integral. The assertion therefore is an immediate consequence
of i. and Prop. 1.1.

\medskip

Recall that $L$ denotes an arbitrary finite extension of $\Qdss_p$
contained in $K$.

\medskip

{\bf Remark 4.6:} {\it $\Zdss_p[[G]], L[[G]]$, and $D_{1/p}(G,K)$
are (left and right) noetherian integral domains with
multiplicative norm $\|\ \|_{1/p}$. }

Proof: Since $L$ is a finitely generated $\Zdss_p$-algebra the
case $L[[G]]$ follows from the case $\Zdss_p[[G]]$. For the other
two cases, using again Prop. 1.1 and Lemma 4.3, it suffices to
know that $gr^\cdot_{1/p} \Zdss_p[[G]]$ is a (left and right)
noetherian integral domain. But this is already contained in [Laz]
since it is shown there, as we have recalled above, that
$gr^\cdot_{1/p}
\Zdss_p[[G]] \cong U(gr(G))$.

\medskip

Turning to flatness properties we will call in the following a
unital ring homomorphism $R \longrightarrow R'$ (faithfully) flat
if $R'$ is (faithfully) flat as a left as well as a right
$R$-module.

\medskip

{\bf Proposition 4.7:} {\it Assuming (HYP) the maps
$$
\Zdss_p[[G]] \longrightarrow L[[G]] \longrightarrow D_r(G,K)
$$
are flat for any $1/p \leq r < 1$ such that $r \in p^\Qdss$. }

Proof: The first map is flat as a base extension of the flat map
$\Zdss_p \hookrightarrow L$. To see that the composite map is flat
it suffices by Prop. 1.2 and Thm. 4.5.i (resp. the argument in the
proof of Remark 4.6 in case $r = 1/p$) to check that the
associated graded map for the filtration $F^\cdot_r$ is flat. This
latter map, by Lemma 4.3, is a base extension of the map
$\Fdss_p[\epsilon]
\longrightarrow gr^\cdot K = k[\epsilon_{\rm o},\epsilon_{\rm o}^{-1}]$
which obviously is flat. To finally deduce the flatness of the
second map it remains to observe that due to the existence of the
trace map $L \longrightarrow \Qdss_p$ the tensor product $D_r(G,K)
\otimes_{L[[G]]} M$, for any $L[[G]]$-module $M$, is naturally a
direct summand of $D_r(G,K) \otimes_{\Qdss_p[[G]]} M = D_r(G,K)
\otimes_{\Zdss_p[[G]]} M$.

\medskip

{\bf Lemma 4.8:} {\it Assume (HYP) and let $1/p < r < 1$ be in
$p^\Qdss$; then $D_{<r}(G,K)$ is (left and right) noetherian and the
map $D_r(G,K) \longrightarrow D_{<r}(G,K)$ is flat. }

Proof: Both assertions can be checked after a faithfully flat base
extension. In particular we may replace the field $K$ by any of its
finite extensions. This allows us to assume that the absolute
ramification index $e$ of $K$ has the property that
$$
r^{\omega(h_i)} = p^{-m_i/e}\ \ \ \hbox{for}\ 1 \leq i \leq d\
\hbox{and appropriate}\ m_i \in \Ndss\ .
$$
As before let $k$ denote the residue class field of $K$ and
$\epsilon_{\rm o}$ the principal symbol of a prime element $\pi$ of
$K$. By Prop. 1.1 and 1.2 it suffices to show that the induced map
between the graded rings $gr^\cdot_r D_{r}(G,K) \longrightarrow
gr^\cdot_r D_{<r}(G,K)$ is a flat map between noetherian rings. The
left hand graded ring has been computed in Thm. 4.5.i. Our above
assumption on $K$ means that the values of the norm $\|\ \|_r$ on
$D_{<r}(G,K)$ lie in $|\pi|^\Zdss \cup \{0\}$. Hence
$$
gr^\cdot_r D_{<r}(G,K) = gr^\cdot K \otimes_k gr^0_r D_{<r}(G,K)
$$
Since $gr^\cdot K$ is a finitely generated $k$-algebra we are
reduced to showing that the map $gr^0_r D_{r}(G,K) \longrightarrow
gr^0_r D_{<r}(G,K)$ is a flat map between noetherian rings. By
definition $F^0_r D_{<r}(G,K)$ consists of all series of
$\sum_\alpha d_\alpha\bfb^\alpha$ where $|d_\alpha|r^{\tau\alpha}
\leq 1$. But $r^{\tau\alpha} = |\pi|^{m_1\alpha_1 + \ldots +
m_d\alpha_d}$. Hence, equivalently, $F^0_r D_{<r}(G,K)$ is all
series $\sum_\alpha e_\alpha
(b_1/\pi^{m_1})^{\alpha_1}\ldots(b_d/\pi^{m_d})^{\alpha_d}$ where
$|e_\alpha| \leq 1$. Mapping this series to\break $\sum_\alpha
(e_\alpha {\rm mod}\pi) u_1^{\alpha_1}\ldots u_d^{\alpha_d}$
induces a bijection
$$
gr^0_r D_{<r}(G,K) \mathop{\longrightarrow}\limits^{\cong}
k[[u_1,\ldots,u_d]]
$$
into the formal power series over $k$ in the variables
$u_1\ldots,u_d$. By Lemma 4.4 this is an isomorphism of rings. Clearly
$gr^0_r D_{r}(G,K)$ corresponds to the polynomial ring
$k[u_1,\ldots,u_d]$. It is well known that the inclusion of the
polynomial ring into the formal power series ring over a field is a
flat map between noetherian rings.

\medskip

{\bf Theorem 4.9:} {\it Assuming (HYP) the map $D_r(G,K)
\longrightarrow D_{r'}(G,K)$ is flat for any
$1/p < r' \leq r < 1$ in $p^\Qdss$. }

Proof: We certainly may assume that $r' < r$. Then the map in question
can be viewed as the composite
$$
D_r(G,K) \longrightarrow D_{<r}(G,K) \longrightarrow D_{r'}(G,K)\ .
$$
Because of Lemma 4.8 it remains to be seen that the second map is
flat. Again we may enlarge the field $K$ if convenient. This time we
assume that the absolute ramification index $e$ of $K$ satisfies
$$
r^{\omega(h_i)} = p^{-m_i/e}\ \ \ \hbox{and}\ \ \ (r')^{\omega(h_i)} =
p^{-m'_i/e} \ \ \ \hbox{for}\ 1 \leq i \leq d\
$$
and appropriate $m_i < m'_i$ in $\Ndss$. Since $D_{<r}(G,K) = \Qdss_p
\otimes_{\Zdss_p} F^0_r D_{<r}(G,K)$ it suffices to show that the
inclusion $F^0_r D_{<r}(G,K) \longrightarrow D_{r'}(G,K)$ is flat.
The advantage of this is that the left term is closed (in fact
c-compact) in the Banach space on the right. This means that on
both sides the filtration $F^\cdot_{r'}$ is complete. By Prop. 1.2
we therefore are reduced to showing that the corresponding
injective graded map $gr^\cdot_{r'} F^0_r D_{<r}(G,K)
\longrightarrow gr^\cdot_{r'} D_{r'}(G,K)$ is a flat map between
noetherian rings. By the assumption on $K$ on both sides a
homogeneous component can be nonzero only if its degree lies in
${1 \over e}\Zdss$. Fix a prime element $\pi$ of $K$, and let $m
\in \Zdss$. Then $F^{m/e}_{r'} F^0_r D_{<r}(G,K)$ consists of all
series $\sum_\alpha e_\alpha
(b_1/\pi^{m'_1})^{\alpha_1}\ldots(b_d/\pi^{m'_d})^{\alpha_d}$
where $|e_\alpha| \leq {\rm min}(|\pi|^m,|\pi|^{\sum_i(m'_i -
m_i)\alpha_i})$. On the other hand $F^{m/e}_{r'} D_{r'}(G,K)$
consists of all such series for which $|e_\alpha| \leq |\pi|^m$
and $|e_\alpha|
\longrightarrow 0$ as $|\alpha| \longrightarrow \infty$. Denoting, as
usual, by $k$ the residue class field of $K$ and by $\epsilon_{\rm o}$
the principal symbol of $\pi$ it follows that
$$
gr^{m/e}_{r'} D_{r'}(G,K) = \epsilon_{\rm o}^m \cdot k[u_1,\ldots,u_d]
\ \ \ \hbox{with}\ u_i := \epsilon_{\rm
o}^{-m'_i}\cdot\sigma_{r'}(b_i)
$$
and
$$
\matrix{
gr^{m/e}_{r'} F^0_r D_{<r}(G,K) =
    \bigoplus \{\epsilon_{\rm o}^m \cdot
    ku_1^{\alpha_1}\ldots u_d^{\alpha_d} : \sum_i (m'_i -
    m_i)\alpha_i \leq m \} =
    \hfill\cr\cr
\bigoplus \{k \epsilon_{\rm o}^{m - \sum_i (m'_i -
    m_i)\alpha_i} (\epsilon_{\rm o}^{m'_1 - m_1}u_1 )^{\alpha_1}
    \ldots (\epsilon_{\rm o}^{m'_d - m_d}u_d )^{\alpha_d} : \sum_i (m'_i -
    m_i)\alpha_i \leq m \}. }
$$
Note that always $\sum_i (m'_i - m_i)\alpha_i \geq 0$. We obtain
$$
gr^{\cdot}_{r'} D_{r'}(G,K) = k[\epsilon_{\rm o},\epsilon_{\rm
o}^{-1},u_1,\ldots,u_d]
$$
(compare Thm. 4.5.i) and
$$
gr^{\cdot}_{r'} F^0_r D_{<r}(G,K) = k[\epsilon_{\rm o},\epsilon_{\rm
o}^{m'_1 - m_1}u_1,\ldots,\epsilon_{\rm o}^{m'_d - m_d}u_d]\ .
$$
Obviously both rings are noetherian and the inclusion between them,
being localization in $\epsilon_{\rm o}$, is flat.

\medskip

{\bf Theorem 4.10:} {\it For any compact locally
$\Qdss_p$-analytic group $G$ which carries a $p$-valuation
satisfying (HYP) the algebra $D(G,K)$ is a Fr\'echet-Stein
algebra. }

Proof: This is a consequence of Thm. 4.5.ii and Thm. 4.9.

\medskip

To construct our sequence of norms which exhibits $D(G,K)$ as a
Fr\'echet-Stein algebra we made the two choices of a $p$-valuation
$\omega$ on $G$ and of a corresponding ordered basis
$(h_1,\ldots,h_d)$ of $G$. Consider another ordered basis
$(h'_1,\ldots,h'_d)$ with respect to the same $\omega$. Put $b'_i :=
h'_i - 1$, $\bfb'^{\alpha} := b_1'^{\alpha_1}b_2'^{\alpha_2}\cdots
b_{d}'^{\alpha_d}$, and $\tau'\alpha := \sum_i \alpha_i\omega(h'_i)$
for $\alpha = (\alpha_1,\ldots,\alpha_d) \in \Ndss_0^d$. Consider the
expansions
$$
\bfb'^{\beta} = \sum_\alpha c_{\beta,\alpha}\bfb^\alpha
$$
in $\Zdss_p[[G]]$. Since the norm $\|\ \|_{1/p}$ on $\Zdss_p[[G]]$
does not depend on the choice of the ordered basis we have
$$
\omega_p(c_{\beta,\alpha}) \geq {\rm max}(0,\tau'\beta - \tau\alpha)\
.
$$
Similarly as in the proof of Lemma 4.1.ii this implies that
$$
|c_{\beta,\alpha}|r^{\tau\alpha} \leq r^{\tau'\beta}
$$
for any $1/p \leq r \leq 1$. Expanding an element $\lambda \in D(G,K)$
with respect to both bases as
$$
\lambda = \sum_\alpha d_\alpha\bfb^\alpha = \sum_\beta
d'_\beta\bfb'^\beta
$$
we have $d_\alpha = \sum_\beta d'_\beta c_{\beta,\alpha}$. It follows
that
$$
\mathop{\rm sup}\limits_{\alpha} |d_\alpha|r^{\tau\alpha} \leq
\mathop{\rm sup}\limits_{\beta} |d'_\beta|r^{\tau'\beta}\ .
$$
For reasons of symmetry we must in fact have equality. This shows that
the norms $\|\ \|_r$ on $D(G,K)$, for $1/p \leq r < 1$, are
independent of the choice of the ordered basis.

\medskip

We complete the list of flatness properties with the following result
which answers positively a question raised at the end of \S3 in [ST3].

\medskip

{\bf Theorem 4.11:} {\it Assuming (HYP) the map $L[[G]]
\longrightarrow D(G,K)$ is faithfully flat.}

Proof: The other case being completely analogous we give in the
following only the argument for left faithful flatness. We also agree
that in this proof all occuring numbers $r$ satisfy, without further
mentioning, the conditions $1/p < r < 1$ and $r \in p^{\Qdss}$.

For the flatness we let $J \subseteq L[[G]]$ be a left ideal. We
have to show that the canonical map
$$
D(G,K) \otimes_{L[[G]]} J \longrightarrow D(G,K)
$$
is injective. The ring $L[[G]]$ being noetherian by Remark 4.6 we
have that the left hand $D(G,K)$-module is finitely presented and
hence coadmissible. It follows that
$$
\matrix{
D(G,K) \otimes_{L[[G]]} J & =
         \mathop{\lim\limits_{\longleftarrow}}\limits_{r}
         D_r(G,K) \otimes_{D(G,K)} (D(G,K) \otimes_{L[[G]]} J) \cr\cr
      & = \mathop{\lim\limits_{\longleftarrow}}\limits_{r}
         D_r(G,K) \otimes_{L[[G]]} J\ .\hfill }
$$
But $D_r(G,K)$ is flat over $L[[G]]$ according to Thm. 4.7 so that
the map $D_r(G,K) \otimes_{L[[G]]} J \longrightarrow D_r(G,K)$ is
injective. Since the projective limit is left exact this implies
the injectivity we want.

For faithful flatness we have to show that $D(G,K)
\otimes_{L[[G]]} M \neq 0$ for any nonzero $L[[G]]$-module $M$. Let $o$ denote the ring
of integers in $L$ and put $o[[G]] := o \otimes_{\Zdss_p}
\Zdss_p[[G]]$. Since $L[[G]] \otimes_{o[[G]]} M = L \otimes_o M = M$
and hence $D(G,K) \otimes_{L[[G]]} M = D(G,K) \otimes_{o[[G]]} M$
it suffices to show that $D(G,K) \otimes_{o[[G]]} M \neq 0$ for
any nonzero $p$-torsionfree $o[[G]]$-module $M$. By the usual
argument this reduces to the claim that $D(G,K)\cdot J \neq
D(G,K)$ for any proper left ideal $J \subseteq o[[G]]$ such that
$o[[G]]/J$ is $p$-torsionfree. Since $o[[G]]$ is noetherian by
Remark 4.6 and $o[[G]] \longrightarrow D(G,K)$ is already known to
be flat the $D(G,K)$-module $D(G,K)\cdot J = D(G,K)
\otimes_{o[[G]]} J$ is coadmissible which implies that
$$
D(G,K)\cdot J = \mathop{\lim\limits_{\longleftarrow}}\limits_{r}
D_r(G,K) \otimes_{o[[G]]} J =
\mathop{\lim\limits_{\longleftarrow}}\limits_{r} D_r(G,K)\cdot J
$$
the second identity coming from the fact that $D_r(G,K)$ is flat
over $o[[G]]$ by Thm. 4.7. It therefore suffices to find at least
one $r$ such that $D_r(G,K)\cdot J \neq D_r(G,K)$.

To simplify the notation we use in the following the abbreviations
$R := o[[G]]$ and $D_r := D_r(G,K)$. Giving the outer terms of the
exact sequence
$$
0 \longrightarrow D_rJ \longrightarrow D_r \longrightarrow D_r/D_rJ
\longrightarrow 0\ .
$$
the filtration induced by $F^\cdot_r D_r$ we obtain the exact sequence
of associated graded $gr^\cdot_r D_r$-modules
$$
0 \longrightarrow gr^\cdot_r D_rJ \longrightarrow gr^\cdot_r D_r
\longrightarrow gr^\cdot_r (D_r/D_rJ)
\longrightarrow 0\ .
$$
We compare this with the exact sequence of $gr^\cdot_r R$-modules
$$
0 \longrightarrow gr^\cdot_r J \longrightarrow gr^\cdot_r R
\longrightarrow gr^\cdot_r (R/J) \longrightarrow 0
$$
arising from giving the outer terms of the exact sequence
$$
0 \longrightarrow J \longrightarrow R \longrightarrow R/J
\longrightarrow 0
$$
the filtration induced by $F^\cdot_r R$. By Lemma 4.3 the map
$gr^\cdot_r R \longrightarrow gr^\cdot_r D_r$ is flat. Hence we have
the commutative exact diagram:
$$
\xymatrix{
  0 \ar[d] & 0 \ar[d] \\
  gr^\cdot_r D_r \otimes_{gr^\cdot_r R} gr^\cdot_r J
  \ar[d] \ar[r]^(0.6){(1)} & gr^\cdot_r D_rJ \ar[d] \\
  gr^\cdot_r D_r \otimes_{gr^\cdot_r R} gr^\cdot_r R
  \ar[d] \ar[r]^(0.6){\cong} & gr^\cdot_r D_r \ar[d] \\
  gr^\cdot_r D_r \otimes_{gr^\cdot_r R} gr^\cdot_r
  (R/J) \ar[d] \ar[r]^(0.6){(2)} & gr^\cdot_r (D_r/D_rJ) \ar[d] \\
  0 & 0 }
$$
The horizontal arrow in the middle obviously is an isomorphism. Since
both rings, $gr^\cdot_r R$ and $gr^\cdot_r D_r$, are noetherian by
Thm. 4.5.i we may apply Lemma 1.3 and obtain that $(1)$ and hence
$(2)$ is an isomorphism. This reduces us to finding an $r$ such that
$gr^\cdot_r D_r \otimes_{gr^\cdot_r R} gr^\cdot_r (R/J) \neq 0$. Since
$R$ is compact the finitely generated left ideal $J$ is closed. Hence
the filtration $F^\cdot_r (R/J)$ is separated and $gr^\cdot_r (R/J)
\neq 0$ if $R/J \neq 0$. From Lemma 4.3 we know that $gr^\cdot_r D_r$
is faithfully flat over $(gr^\cdot_r R)[\epsilon^{-1}]$. So we are
finally reduced to proving that
$$
gr^\cdot_r (R/J)\ \hbox{is}\ \epsilon\ \hbox{-torsionfree for}\ r\
\hbox{sufficiently close to}\ 1
$$
provided $R/J$ is $p$-torsionfree.

Let $\pi \in o$ be a prime element of $L$ and let $e$ denote the
absolute ramification index of $L$. It is more convenient to study
the multiplication by $\sigma_r(\pi)$ (instead of $\epsilon =
\sigma_r(p)$) on $gr^\cdot_r (R/J)$. Suppose therefore that
$$
x \in F^s_r R + J \ \ \hbox{and}\ \ \pi x \in F^{(s + 1/e)+}_r R + J\
.
$$
We have to show that $x \in F^{s+}_r R + J$ (provided $r$ is
sufficiently close to $1$). Write
$$
\pi x = y + a\ \ \ \hbox{with}\ \ y \in F^{(s + 1/e)+}_r R\ \ \hbox{and}\ \
a \in J\ .
$$
We will find $y' \in R$ and $a' \in J \cap F^{(s + 1/e)+}_r R$ such
that $y = \pi y' + a'$. Then $\pi y' = y - a' \in F^{(s + 1/e)+}_r R$,
hence $y' \in F^{s+}_r R$, and $\pi (x - y') = a + a' \in J$. Since
$R/J$ is assumed to be $\pi$-torsionfree it follows that
$$
x \in y' + J \subseteq F^{s+}_r R + J\ .
$$
To actually find $y'$ and $a'$ we look at the left ideal $J + \pi
R/\pi R$ in the ring $R/\pi R$. On $R/\pi R$ we consider the quotient
filtration $F^s_r(R/\pi R) := F^s_r R + \pi R/\pi R$ and then on $J +
\pi R/\pi R$ the induced filtration
$$
F^s_r(J + \pi R/\pi R) := [(J + \pi R) \cap (F^s_r R + \pi R)]/\pi R\
.
$$
Since $R$ is compact the quotient filtration on $R/\pi R$ is complete,
and $gr^\cdot_r(R/\pi R)$ being a quotient of $gr^\cdot_r R$ is
noetherian. Hence we may apply Lemma 1.4 to the ideal $J + \pi R/\pi
R$ in $R/\pi R$. We obtain elements $a_1,\ldots,a_l \in J\setminus \pi
R$ such that $a_i + \pi R \in F^{s_i}_r(R/\pi R) \setminus
F^{s_i+}_r(R/\pi R)$ and
$$
F^s_r(J + \pi R/\pi R) = \sum_{i=1}^l F^{s-s_i}_r (R/\pi R)\cdot (a_i
+
\pi R)
$$
for any $s \in \Rdss$. Since
$$
y + \pi R = - a + \pi R \in F^{(s+1/e)+}_r (J + \pi R/\pi R)
$$
we find $b_i \in F^{(s+1/e - s_i)+}_r R$ such that
$$
y + \pi R = b_1a_1 + \ldots + b_la_l + \pi R
$$
or equivalently
$$
y = \pi y' + a'\ \ \ \hbox{for}\ a' := b_1a_1 + \ldots + b_la_l
\in J\
\hbox{and some}\ y' \in R\ .
$$
The element $a'$ will lie in $F^{(s+1/e)+}_r R$ provided we have $a_i
\in F^{s_i}_r R$ for any $1 \leq i \leq l$. But the subsequent lemma
says that this automatically holds true for any $r$ sufficiently close
to $1$.

\medskip

{\bf Lemma 4.12:} {\it Let $o$ denote the ring of integers in $L$
and $\pi \in o$ a prime element and put $R := o
\otimes_{\Zdss_p} \Zdss_p[[G]]$. For any $a \in R \setminus \pi R$ there is a $1/p \leq
r(a) < 1$ such that for all $r(a) \leq r < 1$ and all $s \in \Rdss$ we
have: If $a \in F^s_r R + \pi R$ then $a \in F^s_r R$. }

Proof: Write $ a = \sum_\alpha d_\alpha\bfb^\alpha$ with $d_\alpha \in
o$. By assumption we have $|d_\alpha| = 1$ for some $\alpha$. We
choose a $\beta$ such that $\tau\beta$ is minimal in the set
$\{\tau\alpha : |d_\alpha| = 1\}$. For any $\alpha$ with $\tau\alpha
\geq \tau\beta$ we have $|d_\alpha|r^{\tau\alpha} \leq r^{\tau\beta}$
for any $r < 1$. Consider on the other hand the finitely many
$\alpha_1,\ldots,\alpha_l$ such that $\tau\alpha_i < \tau\beta$. Since
$|d_{\alpha_i}| < 1$ we find a $1/p \leq r(a) < 1$ such that
$|d_{\alpha_i}| \leq r(a)^{\tau\beta - \tau\alpha_i}$ for any $1 \leq
i \leq l$. It follows that $|d_\alpha|r^{\tau\alpha} \leq
r^{\tau\beta}$ for any $\alpha$ and any $r(a) \leq r < 1$ and hence
that
$$
\|a\|_r \leq r^{\tau\beta}\ \ \ \hbox{for any}\ r(a) \leq r < 1\ .
$$
Suppose now that $a \in F^s_r R + \pi R$ for some $r(a) \leq r < 1$.
Then $|d_\alpha|r^{\tau\alpha} \leq p^{-s}$ for any $\alpha$ such that
$|d_\alpha| = 1$. In particular $r^{\tau\beta} =
|d_\beta|r^{\tau\beta} \leq p^{-s}$ and hence $\|a\|_r \leq p^{-s}$.

\medskip

{\bf 5. Distribution algebras of compact Lie groups}

\smallskip

As before let $\Qdss_p \subseteq L \subseteq K \subseteq \Cdss_p$ be
complete intermediate fields where $L/\Qdss_p$ is finite and $K$ is
discretely valued, and let $G$ be an arbitrary compact locally
$L$-analytic group. We denote by $G_0$ the same group but viewed as a
locally $\Qdss_p$-analytic group.

\medskip

{\bf Theorem 5.1:} {\it $D(G,K)$ is a Fr\'echet-Stein algebra.}

Proof: {\it Step 1:} We first show that $D(G_0,K)$ is a
Fr\'echet-Stein algebra. For this we choose an open normal subgroup
$H_0 \subseteq G_0$ which is a uniform pro-$p$-group. According to
[DDMS] \S4.2 the lower $p$-series (loc. cit. Def. 1.15 and Cor. 1.20)
in $H_0$ is the filtration corresponding to an integrally valued
$p$-valuation $\omega$ on $H_0$ (in fact, for $p = 2$ we have to
replace $H_0$ by its Frattini subgroup which is the first step in its
lower $p$-series so that the restriction of $\omega$ satisfies the
axioms of a $p$-valuation). Let $(h_1,\ldots,h_d)$ be an ordered basis
of $H_0$ corresponding to this $\omega$. Since the lower $p$-series
consists of characteristic subgroups we have
$$
\omega(ghg^{-1}) = \omega(h)\ \ \ \ \hbox{for any}\ g \in G_0, h \in
H_0\ .
$$
Hence $(gh_1g^{-1},\ldots,ghg^{-1})$, for any $g \in G_0$, is an
ordered basis as well. Let $\|\ \|_r$, for $1/p < r < 1$, be the
family of submultiplicative norms on $D(H_0,K)$, corresponding to
$\omega$ and $(h_1,\ldots,h_d)$, which we have constructed in section
4. As discussed after Thm. 4.10 each of these norms is in fact
independent of the choice of the ordered basis. This means that we
have
$$
\|\delta_g \lambda \delta_{g^{-1}}\|_r = \|\lambda\|_r\ \ \ \ \hbox{for
any}\ g \in G_0\ \hbox{and}\ \lambda \in D(H_0,K)\ .
$$
We now fix coset representatives $1 = g_1,g_2,\ldots,g_l$ for
$H_0$ in $G_0$. Then the Dirac distributions
$\delta_{g_1},\ldots,\delta_{g_l}$ form a basis of $D(G_0,K)$ as a
left $D(H_0,K)$-module. For any $\mu
\in D(G_0,K)$ write $\mu = \lambda_1(\mu)\delta_{g_1} + \ldots +
\lambda_l(\mu)\delta_{g_l}$ and put
$$
q_r(\mu) := {\rm
max}(\|\lambda_1(\mu)\|_r,\ldots,\|\lambda_l(\mu)\|_r)\ .
$$
Since $D(G_0,K)$ is topologically isomorphic to $D(H_0,K)^l$ by the
open mapping theorem this defines a continuous norm $q_r$ on
$D(G_0,K)$ which extends the norm $\|\ \|_r$ on $D(H_0,K)$. It also is
clear that the family $q_r$, for $1/p < r < 1$, defines the Fr\'echet
topology of $D(G_0,K)$.

We claim that the multiplication in $D(G_0,K)$ is continuous with
respect to each norm $q_r$. For any two $\mu,\mu' \in D(G_0,K)$ we
compute
$$
\matrix{
q_r(\mu\mu') & = q_r(\mathop{\sum}\limits_{i,j}
\lambda_i(\mu)\delta_{g_i}\lambda_j(\mu')\delta_{g_j})\hfill\cr\cr
& = q_r(\mathop{\sum}\limits_{i,j}
\lambda_i(\mu)(\delta_{g_i}\lambda_j(\mu')
\delta_{g_i^{-1}})\delta_{g_i g_j})\hfill\cr\cr
& = q_r(\mathop{\sum}\limits_k \mathop{\sum}\limits_{i,j}
\lambda_i(\mu)(\delta_{g_i}\lambda_j(\mu')
\delta_{g_i^{-1}})\lambda_k(\delta_{g_i g_j})\delta_{g_k})\hfill\cr\cr
& \leq \mathop{\rm max}\limits_{i,j,k}
\|\lambda_i(\mu)\|_r\|\delta_{g_i}\lambda_j(\mu')
\delta_{g_i^{-1}}\|_r\|\lambda_k(\delta_{g_i g_j})\|_r\hfill\cr\cr
& = \mathop{\rm max}\limits_{i,j,k}
\|\lambda_i(\mu)\|_r\|\lambda_j(\mu')\|_r
\|\lambda_k(\delta_{g_i g_j})\|_r\hfill\cr\cr
& \leq \mathop{\rm max}\limits_{i,j,k} \|\lambda_k(\delta_{g_i
g_j})\|_r \cdot q_r(\mu) \cdot q_r(\mu')\ .\hfill }
$$
It follows that the completion $D_r(G_0,K)$ of $D(G_0,K)$ with
respect to $q_r$ is a $K$-Banach algebra which contains
$D_r(H_0,K)$ and which as a left $D_r(H_0,K)$-module is free with
basis $\delta_{g_1},\ldots,\delta_{g_l}$. Hence Thm. 4.5.ii
implies that $D_r(G_0,K)$ is left noetherian provided $r \in
p^{\Qdss}$. Moreover, since
$$
D_{r'}(G_0,K) = D_{r'}(H_0,K) \mathop{\otimes}\limits_{D_r(H_0,K)}
D_r(G_0,K)
$$
as bimodules we obtain from Thm. 4.9 that $D_{r'}(G_0,K)$ is flat as a
right $D_r(G_0,K)$-module for any $1/p < r' \leq r < 1$ in
$p^{\Qdss}$. (We note that the above argument is left-right symmetric;
in particular $D_r(G_0,K)$ is left and right noetherian if $r \in
p^{\Qdss}$.)

{\it Step 2:} To derive our assertion from Step 1 we consider the
obvious injective map between vector spaces of locally analytic
functions
$$
C^{an}(G,K) \hookrightarrow C^{an}(G_0,K)\ .
$$
It is in fact a topological embedding (see the proof of [ST4] Lemma
1.2). By the Hahn-Banach theorem and the open mapping theorem we
therefore obtain dually a quotient map of Fr\'echet algebras
$$
D(G_0,K) \dlongrightarrow D(G,K)\ .
$$
Hence the assertion follows from the first step and Prop. 3.7.

\medskip

{\bf Theorem 5.2:} {\it The map $L[[G]] \longrightarrow D(G_0,K)$ is
faithfully flat. }

Proof: Choose an open normal subgroup $H_0 \subseteq G_0$ which is a
uniform pro-$p$-group. We then have the bimodule isomorphism
$$
D(G_0,K) \cong L[[G]] \mathop{\otimes}\limits_{L[[H_0]]} D(H_0,K)\ .
$$
Hence the assertion follows by base change (compare [B-CA] I\S3.3
Prop. 5) from Thm. 4.11.

\medskip

It is likely that also the map $L[[G]] \longrightarrow D(G,K)$ is
faithfully flat. But this seems to require additional arguments beyond
the present methods.

\medskip

{\bf 6. Admissible representations}

\smallskip

With $\Qdss_p \subseteq L \subseteq K \subseteq \Cdss_p$ as in the
previous sections we now let $G$ be an arbitrary locally $L$-analytic
group.

\medskip

{\bf Definition:} {\it A (left) $D(G,K)$-module is called coadmissible
if it is coadmissible as a $D(H,K)$-module for every compact open
subgroup $H \subseteq G$.}

\medskip

For any two compact open subgroups $H \subseteq H' \subseteq G$ the
algebra $D(H',G)$ is finitely generated free and hence coadmissible as
a $D(H,K)$-module. It therefore follows from Lemma 3.8 that the
condition in the above definition needs to be tested only for a single
compact open subgroup $H$. We let $\Cscr_G$ denote the full
subcategory of all coadmissible $D(G,K)$-modules in ${\rm
Mod}(D(G,K))$. As a consequence of Cor. 3.4 it is closed with respect
to the formation of finite direct sums, kernels, and cokernels. In
particular $\Cscr_G$ is an abelian category. A module $M$ in $\Cscr_G$
viewed as a $D(H,K)$-module for some compact open subgroup $H
\subseteq G$ carries its canonical Fr\'echet topology. The arguments
in the proof of Lemma 3.8 applied to the algebra homomorphism $D(H,K)
\longrightarrow D(H',K)$ show that this topology in fact is
independent of the choice of $H$. We therefore call it the {\it
canonical topology} on $M$. Consider, for any $g \in G$, the map $M
\mathop{\longrightarrow}\limits^{g\cdot} M$. If we view the source as
a $D(H,K)$-module and the target as a $D(gHg^{-1},K)$-module then this
map is compatible with these module structures relative to the
isomorphism of Fr\'echet algebras $\delta_g . \delta_{g^{-1}} : D(H,K)
\mathop{\longrightarrow}\limits^{\cong} D(gHg^{-1},K)$. Using Lemma
3.8 we conclude that this map $M
\mathop{\longrightarrow}\limits^{g\cdot} M$ is continuous for the
canonical topology. Since
$$
D(G,K) = \bigoplus_{g\in H\setminus G} D(H,K)\delta_g
$$
as locally convex vector spaces it then also follows that the
$D(G,K)$-action on $M$ is separately continuous. Moreover, any map in
$\Cscr_G$ is continuous and strict for the canonical topologies.

\medskip

{\bf Lemma 6.1:} {\it The canonical topology on a coadmissible
$D(G,K)$-module $M$ is nuclear.}

Proof: We choose a compact open subgroup $H \subseteq G$ such that
the locally $\Qdss_p$-analytic group $H_0$ which underlies $H$
carries a $p$-valuation $\omega$ satisfying the condition (HYP) in
section 4. We also fix an ordered basis $(h_1,\ldots,h_d)$ of
$(H_0,\omega)$ and denote by $\|\ \|_r$, for $1/p \leq r < 1$, the
family of norms on $D(H_0,K)$ constructed in section 4. The
completion of $D(H_0,K)$, resp. of the quotient algebra $D(H,K)$,
with respect to $\|\ \|_r$, resp. the quotient norm of $\|\ \|_r$,
is denoted by $D_r(H_0,K)$, resp. $D_r(H,K)$. As a Fr\'echet space
$M$ is the projective limit
$$
M = \mathop{\lim\limits_{\longleftarrow}}\limits_{r} M_r
$$
of the Banach spaces
$$
M_r := D_r(H,K) \mathop{\otimes}\limits_{D(H,K)} M\ .
$$
For its nuclearity it suffices to check (compare [NFA] Prop. 19.9)
that the transition maps $M_r \longrightarrow M_{r'}$, for $1/p
\leq r' < r < 1$, are compact. But we have, for some $l \in
\Ndss$, a commutative diagram
$$
\xymatrix{
  D_r(H_0,K)^l \ar[d] \ar@{->>}[r] & D_r(H,K)^l \ar[d] \ar@{->>}[r] & M_r \ar[d] \\
  D_{r'}(H_0,K)^l \ar@{->>}[r] & D_{r'}(H,K)^l \ar@{->>}[r] & M_{r'}   }
$$
where the horizontal maps are strict surjections. Hence it
suffices to show that the first perpendicular arrow is compact. In
other words we are reduced to showing that the natural map
$$
D_r(H_0,K) \longrightarrow D_{r'}(H_0,K)
$$
is compact. The ring structure being irrelevant here we may use
the global chart $\psi$ corresponding to the fixed ordered basis
to further reduce this to the compactness of the map
$$
D_r(\Zdss_p^d,K) \longrightarrow D_{r'}(\Zdss_p^d,K)\ .
$$
As we have discussed after Prop. 4.2, both sides are Banach spaces
of convergent power series on certain closed disks strictly
contained in each other and the map is the restriction map. By a
well known nonarchimedean version of Montel's theorem (compare
[NFA] p.99) this map is compact.

\medskip

{\bf Remark 6.2:} {\it Suppose that $G$ is compact; then the analytic
$D(G,K)$-modules in the sense of [ST1] \S1 are precisely the finitely
generated coadmissible $D(G,K)$-modules.}

Proof: If $M$ is a finitely generated and coadmissible $D(G,K)$-module
then, in view of its canonical Fr\'echet topology, it follows from
[ST1] Prop. 1.1 that $M$ is analytic. Vice versa, if $M$ is analytic
then it is shown in the proof of loc. cit. that $M$ is isomorphic to a
quotient of a finitely generated free $D(G,K)$-module (with its
canonical topology) by a closed submodule. Hence Lemma 3.6 implies
that $M$ is coadmissible.

\medskip

It follows from the above remark (and Cor. 3.4.v) that for compact $G$
any finitely presented $D(G,K)$-module is analytic. This was posed as
an open question in [ST1].

We recall from [ST2] \S3 that a locally analytic
$G$-representation $V$ (over $K$) is a barrelled locally convex
Hausdorff $K$-vector space $V$ equipped with a $G$-action by
continuous linear endomorphisms such that, for each $v\in V$, the
orbit map $\rho_{v}(g):=gv$ is a $V$-valued locally analytic
function on $G$. The $G$-action extends to a separately continuous
action of the algebra $D(G,K)$ on $V$. It was proved in [ST2] Cor.
3.3 that the functor $V \longmapsto V'_b$ of passing to the strong
dual induces an anti-equivalence of categories:
$$
\matrix{
\matrix{
\hbox{\rm locally analytic $G$-representations}\cr
\hbox{\rm on $K$-vector spaces of compact type}\cr
\hbox{\rm with continuous linear $G$-maps}\cr
\cr
} &\to&
\matrix{
\hbox{\rm separately continuous $D(G,K)$-}\cr
\hbox{\rm modules on nuclear Fr\'echet}\cr
\hbox{\rm spaces with continuous}\cr
\hbox{\rm $D(G,K)$-module maps}\cr
}}
$$
The objects on the right hand side naturally come as right
$D(G,K)$-modules. But using the anti-involution $g \longmapsto
g^{-1}$ we may and will view them as left $D(G,K)$-modules.

\medskip

{\bf Definition:} {\it An admissible $G$-representation over $K$ is a
locally analytic $G$-representation on a $K$-vector space of compact
type $V$ such that the strong dual $V'_b$ is a coadmissible
$D(G,K)$-module equipped with its canonical topology.}

\medskip

We let ${\rm Rep}_K^a(G)$ denote the category of all admissible
$G$-representations over $K$ with continuous $K$-linear $G$-maps.

\medskip

{\bf Theorem 6.3:} {\it The functor
$$
\matrix{
{\rm Rep}_K^a(G) & \mathop{\longrightarrow}\limits^{\sim} & \Cscr_G\cr
\hfill V & \longmapsto & V'\hfill }
$$
is an anti-equivalence of categories. }

Proof: Consider the anti-equivalence which we have recalled above from
[ST2]. According to Lemma 6.1 equipping a coadmissible $D(G,K)$-module
with its canonical topology makes $\Cscr_G$ into a full subcategory of
the right hand side which, by definition, is the essential image of
${\rm Rep}_K^a(G)$ under the functor $V \longmapsto V'_b$.

\medskip

{\bf Proposition 6.4:} {\it i. ${\rm Rep}_K^a(G)$ is an abelian
category; kernel and image of a morphism in ${\rm Rep}_K^a(G)$ are
the algebraic kernel and image with the subspace topology;

ii. any map in ${\rm Rep}_K^a(G)$ is strict and has closed image;

iii. the category ${\rm Rep}_K^a(G)$ is closed with respect to the
passage to closed $G$-invariant subspaces.}

Proof: ${\rm Rep}_K^a(G)$ is abelian as a consequence of Thm. 6.3.
Let $f : V \longrightarrow W$ be a map in ${\rm Rep}_K^a(G)$. The
dual map $f' : W'_b \longrightarrow V'_b$, as a map between
coadmissible modules with their canonical topology, is a strict
map between nuclear and hence ([NFA] Cor. 19.3.ii) reflexive
Fr\'echet spaces. It then follows from the nonarchimedean analog
of [B-TVS] IV\S4.2 Cor. 3 and the fact that over a nonarchimedean
field every reflexive space already is Montel ([NFA] Prop. 15.3)
that $f = f''$ is strict as well and has closed image. The
statement about the kernel and image of $f$ is an easy consequence
of this. The assertion iii. finally follows from Lemma 3.6 and
[ST2] Prop. 1.2(i) (compare the proof of [ST2] Lemma 3.5).

\medskip

In [ST2] \S3 we have introduced, for compact groups $G$, the notion of
a strongly admissible $G$-representation. This is a locally analytic
$G$-representation in a vector space of compact type $V$ such that
$V'$ is finitely generated as a $D(G,K)$-module. But then $V'_b$ is a
quotient of some $D(G,K)^l$ by some closed submodule and hence is
coadmissible (with its canonical topology) by Lemma 3.6. It therefore
follows from Thm. 6.3 that the strongly admissible $G$-representations
are precisely the admissible $G$-representations for which $V'$ is a
finitely generated $D(G,K)$-module (provided $G$ is compact).

The following technical useful observation says that, for compact
$G$ again, the structure of a vector space of compact type which
underlies an admissible $G$-representation $V$ in fact can be
chosen in a $G$-equivariant way.

\medskip

{\bf Proposition 6.5:} {\it Suppose that $G$ is compact; any
admissible $G$-representation $V$ is a compact inductive limit of
locally analytic $G$-representations on $K$-Banach spaces $V_n$
for $n \in \Ndss$.}

Proof: Let $(r_n)_{n\in\Ndss}$ be a sequence of real numbers in
$\{r \in p^{\Qdss} : 1/p < r < 1\}$ which converges to $1$, and
let $A_n := D_{r_n}(G,K)$ be the corresponding sequence of Banach
algebras constructed in the proof of Thm. 5.1 which realizes the
structure of $D(G,K)$ as a Fr\'echet-Stein algebra. The Fr\'echet
space $V'_b$ is the projective limit of the projective system
$$
\ldots \longrightarrow M_n \longrightarrow\ldots\longrightarrow
M_1
$$
of Banach spaces $M_n := A_n \otimes_{D(G,K)} V'$. We let $V_n :=
M'_n$ denote the dual Banach spaces. It was shown in the proof of
Lemma 6.1 that the transition maps in the above projective system
are compact (note that as Banach spaces the $M_n$ remain unchanged
if $G$ is replaced by an open normal subgroup as we do in that
proof). Moreover, by Theorem A, the projection maps $V'
\longrightarrow M_n$ have dense image. It follows that the maps in
the dual inductive system
$$
V_1 \longrightarrow\ldots\longrightarrow V_n \longrightarrow\ldots
$$
are injective and compact ([NFA] Lemma 16.4) and that ([NFA] Prop.
16.5)
$$
V = (V'_b)'_b = \mathop{\lim\limits_{\longrightarrow}}\limits_{n}
V_n
$$ is the compact inductive limit of the $V_n$. By construction
the original projective system is a system of continuous
$D(G,K)$-modules. By functoriality the dual inductive system
therefore is $G$-equivariant. It remains to show that the
$G$-action on each $V_n$ is locally analytic. The continuous
$A_n$-module structure on $M_n$ amounts to a continuous
homomorphism $A_n \longrightarrow \Lscr_s(M_n,M_n)$ into the space
of bounded operators on $M_n$ equipped with the topology of
pointwise convergence. By [B-TVS] III.31 Prop. 6 this map $A_n
\longrightarrow \Lscr_b(M_n,M_n)$ remains continuous when the
target is given the operator norm topology . The dual action of
$A_n$ on $V_n$ arises by composing this map with the passage to
the adjoint and therefore is continuous as well. Hence the
$G$-action on $V_n$ is the restriction of a continuous
$D(G,K)$-action. Moreover, the orbit map $\rho_v$, for any $v \in
V_n$, is the image of the continuous linear map $\lambda
\longmapsto \lambda v$ in $\Lscr(D(G,K),V_n)$ under the
integration isomorphism ([ST2] Thm. 2.2)
$$
\Lscr(D(G,K),V_n) \cong C^{an}(G,V)
$$
and consequently is locally analytic.

\medskip

Next we want to discuss how the smooth representation theory of
$G$ relates to the category ${\rm Rep}_K^a(G)$. We recall that a
smooth $G$-representation $V$ (over $K$) is a $K$-vector space $V$
with a linear $G$-action such that the stabilizer of each vector
in $V$ is open in $G$. A smooth $G$-representation $V$ is called
admissible-smooth if, for any compact open subgroup $H \subseteq
G$, the vector subspace $V^H$ of $H$-invariant vectors in $V$ is
finite dimensional. As we discussed in [ST1] \S2, any
admissible-smooth $G$-representation $V$ equipped with the finest
locally convex topology is a locally analytic $G$-representation
on a vector space of compact type. The derived action of the Lie
algebra $\gfr$ of $G$ in this case obviously is trivial.

\medskip

{\bf Theorem 6.6:} {\it i. Any admissible-smooth
$G$-representation $V$ equipped with the finest locally convex
topology is admissible and has a trivial derived Lie algebra
action;

ii. any admissible $G$-representation $V$ with trivial derived Lie
algebra action is admissible-smooth equipped with the finest
locally convex topology.}

Proof: We fix a compact open subgroup $H \subseteq G$. Let
$I(\gfr) \subseteq D(H,K)$ denote the closed two sided ideal
generated by $\gfr$. The quotient algebra $D^{\infty}(H,K) :=
D(H,K)/I(\gfr)$ is the algebra of locally constant distributions
on $H$; it is the projective limit
$$
D^{\infty}(H,K) = \mathop{\lim\limits_{\longleftarrow}}\limits_{N}
K[H/N]
$$
of the algebraic group rings $K[H/N]$ where $N$ runs over the open
normal subgroups of $H$ (see [ST1] \S2). This latter description
visibly exhibits $D^{\infty}(H,K)$ as a $K$-Fr\'echet-Stein
algebra.

i. The strong dual of $V$ when we equip $V$ with the finest locally
convex topology is the full linear dual $\Hom_K(V,K)$ with the
topology of pointwise convergence. Since $\gfr$ acts trivially on $V$
the $D(H,K)$-action on $\Hom_K(V,K)$ factorizes through
$D^{\infty}(H,K)$. As $V$ is the locally convex inductive limit of the
finite dimensional vector spaces $V^N$ its strong dual satisfies
$$
V'_b = \mathop{\lim\limits_{\longleftarrow}}\limits_{N}
\Hom_K(V^N,K)\ .
$$
In addition we have
$$
K[H/N] \mathop{\otimes}\limits_{D^{\infty}(H,K)} \Hom_K(V,K) =
\Hom_K(V^N,K)\ .
$$
This shows that $V'_b$ is coadmissible with the canonical topology
as a $D^{\infty}(H,K)$- and hence as a $D(H,K)$-module.

ii. The $D(H,K)$-action on $V$ and hence on $V'_b$ factorizes
through $D^{\infty}(H,K)$ by assumption. Hence the $G$-action on
$V$ is smooth ([ST1] \S2). Secondly, as a coadmissible
$D^{\infty}(H,K)$-module with the canonical topology $V'_b$ is the
projective limit of the finite dimensional vector spaces $K[H/N]
\otimes_{D^{\infty}(H,K)} V'$. Dualizing again we obtain that $V
\cong (V'_b)'_b$ is a locally convex inductive limit of finite
dimensional vector spaces and therefore carries the finest locally
convex topology. The finite dimensionality of
$$
K[H/N] \mathop{\otimes}\limits_{D^{\infty}(H,K)} V' = K[H/N]
\mathop{\otimes}\limits_{D^{\infty}(H,K)} \Hom_K(V,K) =
\Hom_K(V^N,K)
$$
finally implies that $V^N$ is finite dimensional for each $N$.
Hence the $G$-action on $V$ is admissible-smooth.

\medskip

A vector $v \in V$ in a general admissible $G$-representation $V$
is called {\it smooth} if the orbit map $\rho_v$ is locally
constant on $G$. The smooth vectors form a $G$-invariant vector
subspace $V_{smooth}$ in $V$. Taylor's formula shows that a vector
$v$ is smooth if and only if $v$ is annihilated by the derived
$\gfr$-action on $V$. Hence $V_{smooth}$ coincides with the
largest vector subspace of $V$ on which the $\gfr$-action is
trivial:
$$
V_{smooth} = V^{\gfr = 0}
$$
In particular, $V_{smooth}$ is closed in $V$ and is therefore an
admissible $G$-representa- tion by Prop. 6.4.iii. Thm. 6.6.ii then
implies that $V_{smooth}$ in fact is the largest admissible-smooth
subrepresentation of $V$.

\medskip

{\bf 7. The existence of analytic vectors}

\smallskip

In this section we will explore the faithful flatness result Thm.
5.2. Throughout we let $G$ be a compact locally $\Qdss_p$-analytic
group and we let $K$ be a finite extension of $\Qdss_p$ and $o$
its ring of integers. A $K$-Banach space representation $V$ of $G$
is a $K$-Banach space $V$ together with a continuous linear action
$G \times V \longrightarrow V$ of $G$ on $V$. The $G$-action
extends naturally to an action of the algebra $K[[G]]$ by
continuous linear endomorphisms on $V$; the induced action $o[[G]]
\times V \longrightarrow V$ of the compact subring $o[[G]]$ is
continuous (see [ST3] \S2). By functoriality $K[[G]]$ then also
acts on the continuous dual $V'$. In [ST3] \S2 (in particular
Lemma 3.4) we have called a $K$-Banach space representation $V$ of
$G$ {\it admissible} if the dual $V'$ is finitely generated as a
$K[[G]]$-module. As in loc. cit. we let ${\rm Ban}_G^{adm}(K)$
denote the category of all admissible $K$-Banach space
representations of $G$ with morphisms being all $G$-equivariant
continuous linear maps. According to [ST3] Thm. 3.5 the functor
$$
\matrix{
{\rm Ban}_G^{adm}(K) & \mathop{\longrightarrow}\limits^{\sim} &
\Mscr_{K[[G]]}\cr \hfill V & \longmapsto & V'\hfill }
$$
is an anti-equivalence of categories.

Since the algebra $K[[G]]$ is noetherian we have the functor
$$
D(G,K) \mathop{\otimes}\limits_{K[[G]]} . : \Mscr_{K[[G]]}
\longrightarrow \Cscr_G\ .
$$
By composition we obtain a natural functor
$$
\matrix{
{\rm Ban}_G^{adm}(K) & \longrightarrow & {\rm Rep}_K^a(G)\hfill\cr
\hfill V & \longmapsto & (D(G,K) \mathop{\otimes}\limits_{K[[G]]}
V')'_b\ ; }
$$
here the strong dual on the right hand side is formed with respect
to the canonical topology on the coadmissible $D(G,K)$-module
$D(G,K) \otimes_{K[[G]]} V'$. The goal of this section is to give
a description of this functor directly in terms of the $K$-Banach
space representation $V$.

Let us fix a $K$-Banach space representation $V$ of $G$.

\medskip

{\bf Definition:} {\it A vector $v \in V$ is called analytic if
the continuous orbit map $\rho_v$ is locally analytic.}

\medskip

We let $V_{an} \subseteq V$ denote the vector subspace of analytic
vectors. It is clearly $G$-invariant. Moreover the $G$-equivariant
linear map (with $G$ acting by left translations on the right hand
side)
$$
\matrix{
V_{an} & \longrightarrow & C^{an}(G,V)\hfill\cr \hfill v &
\longmapsto & [g \mapsto g^{-1}v] }\leqno{(1)}
$$
is injective. We always equip $V_{an}$ with the subspace topology
with respect to this embedding. (Note that this topology in
general is strictly finer than the topology induced by $V$.) Of
course the $G$-action on $V_{an}$ is continuous.

\medskip

{\bf Theorem 7.1:} {\it For any admissible $K$-Banach space
representation $V$ we have:

i. $V_{an}$ is dense in $V$;

ii. $V_{an}$ is a strongly admissible $G$-representation;

iii. $(V_{an})' = D(G,K) \mathop{\otimes}\limits_{K[[G]]} V'$.

Moreover, the functor
$$
\matrix{ {\rm Ban}_G^{adm}(K) & \longrightarrow & {\rm
Rep}_K^a(G)\hfill\cr \hfill V & \longmapsto & V_{an} \hfill }
$$
is exact.}

Proof: ii. Let $C(G,K)$ denote the $K$-Banach space of continuous
functions on $G$. Given finitely many $l_1,\ldots,l_m \in V'$ the
$G$-equivariant continuous linear map
$$
\matrix{
V & \longrightarrow & C(G,K)^m \hfill\cr v &
\longmapsto & ([g \mapsto l_i(g^{-1}v)])_{1 \leq i \leq m}
}\leqno{(2)}
$$
dualizes into the $K[[G]]$-module homomorphism
$$
\matrix{
K[[G]]^m & \longrightarrow & V' \hfill\cr (\mu_1,\ldots,\mu_m) &
\longmapsto & \mu_1l_1 + \ldots + \mu_ml_m\ . }\leqno{(2)'}
$$
Since, by assumption, $V'$ is finitely generated as a
$K[[G]]$-module it follows, using the Hahn-Banach theorem, that we
may choose the $l_1,\ldots,l_m$ in such a way that they are
$K$-linearly independent and that the map (2) is injective. Its
restriction
$$
\matrix{
V_{an} & \longrightarrow & C^{an}(G,K)^m \hfill\cr v &
\longmapsto & ([g \mapsto l_i(g^{-1}v)])_{1 \leq i \leq m}
}\leqno{(2_{an})}
$$
then is injective as well. Let us consider the commutative
triangle
$$
\xymatrix@R=0.5cm{
                &         C^{an}(G,V) \ar[dd]     \\
  V_{an} \ar@{^{(}->}[ur]^{(1)} \ar@{^{(}->}[dr]_{(2_{an})}                 \\
                &         C^{an}(G,K)^m                 }
$$
where the perpendicular arrow is induced by the linear forms
$l_1,\ldots,l_m$. By the definition of the topology on $V_{an}$
and by [ST5]  Remark 3.6 the map $(1)$ is a closed embedding. The
perpendicular map is continuous and, having a continuous section,
closed. It follows that the map $(2_{an})$ is a closed embedding.
But $C^{an}(G,K)^m$, according to [Fea] 2.3.2 (alternatively use
an argument as in the proof of Lemma 6.1), is a vector space of
compact type. Using [Fea] 3.1.7 and [ST2] Prop. 1.2(i) we conclude
that $V_{an}$ is a locally analytic $G$-representation on a vector
space of compact type. Moreover, by the Hahn-Banach theorem again,
the map $(2_{an})$ dualizes into a $D(G,K)$-linear continuous
surjection between Fr\'echet spaces
$$
D(G,K)^m \dlongrightarrow (V_{an})_b'\ .\leqno{(2_{an})'}
$$
In this situation Lemma 3.6 implies that $(V_{an})_b'$ is a finitely
generated coadmissible $D(G,K)$-module with its canonical topology.

iii. Keeping the above notations we have the commutative diagram:
$$
\xymatrix{
  K[[G]]^m \ar[d] \ar[r]^{(2)'}
                & V' \ar[d] \ar@/^1pc/[ddr]  \\
  D(G,K)^m \ar[r] \ar@/_1pc/[drr]_{(2_{an})'}
                & D(G,K) \otimes_{K[[G]]} V' \ar@{-->}[dr]^{(\ast)}  \\
                &               & (V_{an})'              }
$$
The canonical map $(\ast)$ which, by ii., lies in $\Cscr_G$, is
surjective since the map $(2_{an})'$ is surjective. It is continuous
with respect to the canonical topologies on both sides. Note that the
target with the canonical topology coincides with $(V_{an})'_b$. We
define
$$
V_1 := (D(G,K) \mathop{\otimes}\limits_{K[[G]]} V')'_b
$$
as an admissible $G$-representation. Using Hahn-Banach we see that
$(\ast)$ is injective if the dual map $V_{an} \longrightarrow V_1$ is
surjective. We will establish this by looking at the right hand
perpendicular arrow. Let $M \subseteq V'$ denote the image of
$o[[G]]^m$ under the map $(2)'$. The Fr\'echet topology of $D(G,K)$
induces on $o[[G]]$ its natural compact topology as a completed group
ring. The quotient topology from $o[[G]]^m$ on $M$ is the unique
(compact) Hausdorff topology on $M$ which makes the module structure
$o[[G]] \times M \longrightarrow M$ continuous (compare [ST3] Prop.
3.1). We write $V'_{bs}$ for $V'$ equipped with the finest locally
convex topology such that the inclusion $M \subseteq V'$ is
continuous. Then the map under consideration is continuous as a map of
locally convex vector spaces
$$
V'_{bs} \longrightarrow D(G,K) \mathop{\otimes}\limits_{K[[G]]} V'
$$
(where the target, as always, carries its canonical topology). Since
$K[[G]]$ is dense in $D(G,K)$ ([ST2] Lemma 3.1) this map has dense
image. Hence the dual map $V_1 \longrightarrow (V'_{bs})'_b$ is
injective. The discussion before Lemma 1.4 in [ST3] tells us that as
vector spaces we have $(V'_{bs})' = V$. We therefore have the two
injective maps
$$
V_{an} \longrightarrow V_1 \longrightarrow (V'_{bs})' = V
$$
whose composition obviously is the inclusion $V_{an} \subseteq V$. The
left map is continuous by construction. If we show that the right map
is continuous with respect to the original Banach space topology on
$V$ then its image has to fall into $V_{an}$ since $V_1$ is a locally
analytic $G$-representation, and we obtain $V_{an} = V_1$ as claimed.
For this continuity it suffices to check that the strong topology on
$(V'_{bs})'$ is finer than the Banach space topology on $V$. According
to [ST3] the Banach space $V$ is identified with $\Hom_o^{\rm
cont}(M,K)$ with norm given by the sup-norm, and the unit ball is
$\Hom_o^{\rm cont}(M,o)$. But $M$ is compact and therefore bounded in
$V'_{bs}$. Hence $\Hom_o^{\rm cont}(M,o)$ is open in $(V'_{bs})'_b$.

i. In view of what we know already it suffices to show that the image
of this latter continuous map $V_1 \longrightarrow V$ is dense in $V$.
But this is, by Hahn-Banach, a consequence of the injectivity of the
dual map $V' \longrightarrow D(G,K) \otimes_{K[[G]]} V'$ which in turn
follows from the faithful flatness of $D(G,K)$ over $K[[G]]$ (Thm.
5.2).

Finally, in order to establish the last assertion, we now have
that the functor $V \longmapsto V_{an}$ of passing to the subspace
of analytic vectors coincides with the composite functor:
$$
\matrix{
{\rm Ban}_G^{adm}(K) & \longrightarrow & \Mscr_{K[[G]]} &
\longrightarrow & \Cscr_{D(G,K)} & \longrightarrow & {\rm
Rep}_K^a(G)\cr\cr V & \longmapsto & V' & \longmapsto & D(G,K)
\mathop{\otimes}\limits_{K[[G]]} V' & \longmapsto & (D(G,K)
\otimes V')'_b }
$$
In this composition all three functors are exact: The first one by
the strictness of all maps in ${\rm Ban}_G^{adm}(K)$ ([ST3]) and
Hahn-Banach, the second one by Thm. 5.2, and the third one by the
strictness of all maps in $\Cscr_{D(G,K)}$ with respect to the
canonical topologies and Hahn-Banach. Hence the functor $V
\longmapsto V_{an}$ on ${\rm Ban}_G^{adm}(K)$ is exact.

\medskip

The content of this section generalizes in a straightforward way to
arbitrary locally $\Qdss_p$-analytic groups $G$ if one extends the
notion of an admissible Banach representation from the compact to the
general case by a procedure similar to the one we have used in section
6.

\medskip

{\bf 8. Dimension theory for coadmissible modules}

\smallskip

Over commutative noetherian rings the size of a finitely generated
module can be measured by the dimension of its support. If the
ring, in addition, is regular then this dimension function can be
described in terms of the vanishing of certain Ext-modules. In
this latter form dimension theory was generalized to so called
Auslander regular noncommutative noetherian rings (compare [Bjo]).
The Banach algebras $D_r(G,K)$ which we have introduced in section
4 turn out to be Auslander regular. We will demonstrate in the
present section that, as a consequence, the category $\Cscr_G$,
for any locally $\Qdss_p$-analytic group $G$, allows a well
behaved dimension theory. We will proceed in an axiomatic way, and
we begin by setting up the general formalism of Ext-functors in
our context.

We always will assume in this section that $A$ is a ''two sided''
$K$-Fr\'echet-Stein algebra by which we mean the existence of a
sequence $q_1 \leq\ldots\leq q_n \leq\ldots$ of continuous algebra
seminorms on $A$ defining the Fr\'echet topology of $A$ such
that\hfill\break
 (i) $\ \ A_{q_n}$ is left and right noetherian, and\hfill\break
 (ii) $\ A_{q_n}$ is flat as a left and as a right $A_{q_{n+1}}$-module\hfill\break
for any $n\in\Ndss$.

\medskip

{\bf Definition:} {\it A (left module) sheaf for $(A,(q_n))$ is a
family $S = (S_n)_{n\in\Ndss}$ of left unital $A_{q_n}$-modules
$S_n$ together with $A_{q_{n+1}}$-module homomorphisms $S_{n+1}
\longrightarrow S_n$ for any $n \in \Ndss$.}

\medskip

These sheaves, together with the obvious notion of a homomorphism,
form an abelian category $Sh_{(A,(q_n))}$. It contains as a full
abelian subcategory the category $Coh_{(A,(q_n))}$ of coherent
sheaves as introduced in section 3. The left exact global section
functor is defined by
$$
\matrix{ \Gamma : & Sh_{(A,(q_n))} & \longrightarrow & {\rm
Mod}(A) \cr\cr
 & \hfill S & \longmapsto &
 \mathop{\lim\limits_{\longleftarrow}}\limits_{n} S_n\ .\hfill }
$$
According to Cor. 3.3 its restriction to $Coh_{(A,(q_n))}$ is
exact and induces an equivalence of categories
$$
\Gamma : Coh_{(A,(q_n))} \mathop{\longrightarrow}\limits^{\sim}
\,\Cscr_A\ .
$$
There are obvious right module versions of these notions and facts
which we will need later on and for which we will use
corresponding notations decorated with a superscript ''r''. The
following three facts elaborate on the methods in section II of
[Ban].

\medskip

{\bf Lemma 8.1:} {\it i. The category $Sh_{(A,(q_n))}$ has enough
injective objects;

ii. for any $n_0 \in \Ndss$ the functor
$$
\matrix{ Sh_{(A,(q_n))} & \longrightarrow & {\rm Mod}(A_{q_{n_0}})
\cr
 \hfill (S_n)_n & \longmapsto & S_{n_0} \hfill }
$$
respects injective objects;

iii. the global section functor $\Gamma$ respects injective
objects;

iv. given a coherent sheaf $C = (C_n)_n$ the left exact functor
$$
\matrix{ Sh_{(A,(q_n))} & \longrightarrow & \Ndss\hbox{-projective
systems of abelian groups}\hfill\cr\cr (S_n)_n & \longmapsto &
(\Hom_{A_{q_n}}(C_n,S_n))_n = (\Hom_A(\Gamma(C),S_n)_n\ . }
$$
transforms injective objects into projective systems acyclic for
the functor $\mathop{\rm lim}\limits_{\longleftarrow}$; its $l$-th
derived functor is equal to $(\Ext^l_{A_{q_n}}(C_n,.))_n$.}

Proof: i. The abelian category $Sh_{(A,(q_n))}$ clearly satisfies
$AB5$ and $AB3^{\ast}$. Moreover it is easy to verify that the
family of sheaves $\{(A^{(m)}_n)_n\}_{m\in\Ndss}$ with $A^{(m)}_n
:= 0$, resp. $ := A_{q_n}$, for $n > m$, resp. $n \leq m$, is a
family of generators for $Sh_{(A,(q_n))}$. But any abelian
category with these properties has enough injective objects
(compare [Pop] Thm. 3.10.10).

ii. The functor in question has the exact left adjoint functor $N
\longmapsto (N_n)_n$ with $N_n := 0$ for $n > n_0$ and $N :=
A_{q_n} \otimes_{A_{q_{n_0}}} N$ for $n \leq n_0$.

iii. Let $J = (J_n)_n$ be an injective object in $Sh_{(A,(q_n))}$.
We have to show that, given a left ideal $\lfr \subseteq A$, any
$A$-module homomorphism $\lfr \longrightarrow \Gamma(J)$ extends
to an $A$-module homomorphism $A \longrightarrow \Gamma(J)$. It
follows from Remark 3.2 that $L := (A_{q_n} \otimes_A \lfr)_n$ is
a coherent subsheaf of $(A_{q_n})_n$. The homomorphism $\lfr
\longrightarrow \Gamma(J)$ induces homomorphisms $A_{q_n}
\otimes_A \lfr \longrightarrow A_{q_n} \otimes_A \Gamma(J) = J_n$
which constitute a map $L \longrightarrow J$ in $Sh_{(A,(q_n))}$.
By the injectivity of $J$ this latter map extends to a map
$(A_{q_n})_n \longrightarrow J$ which on global sections induces
the required extension $A \longrightarrow \Gamma(J)$.

iv. The second part of the assertion is a consequence of ii. To
establish the first part we introduce, for any $m \in \Ndss$, the
two functors $S \longmapsto S^{(m)}$ and $S \longmapsto S^{[m]}$
from the category $Sh_{(A,(q_n))}$ into itself defined by
$$
S^{(m)}_n := \cases {S_m\,\quad\hbox{\rm if $n \geq m$,}\cr
                     S_n\ \quad\hbox{\rm if $n \leq m$}}
$$
and
$$
S^{[m]}_n := \cases {0\ \ \quad\hbox{\rm if $n > m$,}\cr
                     S_n \quad\hbox{\rm if $n \leq m$,}}
$$
respectively. The second functor is exact and left adjoint to the
first one. Hence the first functor respects injective objects. We
consider now a fixed injective object $J = (J_n)_n$ in
$Sh_{(A,(q_n))}$. Then the direct product $\prod_m J^{(m)}$ in
$Sh_{(A,(q_n))}$ is injective as well. There are obvious sheaf
maps $\sigma_m : J\longrightarrow J^{(m)}$ and $\sigma^m_{m+1} :
J^{(m+1)} \longrightarrow J^{(m)}$. We claim that the sequence
$$
\xymatrix@C=0.5cm{
  0 \ar[r] & J \ar[rrr]^{\prod\sigma_m\ \ \ \ } &&& \prod_m J^{(m)}
  \ar[rrrr]^{\prod ({\rm id}_{J^{(m)}} - \sigma^m_{m+1})}
  &&&& \prod_m J^{(m)} \ar[r] & 0 }
$$
is exact. This amounts to the exactness of the sequences
$$
0 \longrightarrow J_n \longrightarrow \prod_{m < n} J_m \times
\prod_{m \geq n} J_n \longrightarrow \prod_{m < n} J_m \times
\prod_{m \geq n} J_n \longrightarrow 0
$$
for each $n \in \Ndss$ which is equivalent to the obvious fact
that the projective systems
$$
\ldots \mathop{\longrightarrow}\limits^= J_n
\mathop{\longrightarrow}\limits^= \ldots
\mathop{\longrightarrow}\limits^= J_n \longrightarrow J_{n-1}
\longrightarrow \ldots \longrightarrow J_1
$$
have $J_n$ as projective limit and have vanishing
${\mathop{\lim}\limits_{\longleftarrow}}^{(1)}$-term. Applying the
functor $\Hom_{Sh_{(A,(q_n))}}(C,.)$ to the above short exact
sequence of injective sheaves we obtain the short exact sequence
$$
0 \rightarrow \Hom_{Sh_{(A,(q_n))}}(C,J) \rightarrow \prod_m
\Hom_{A_{q_m}}(C_m,J_m) \rightarrow \prod_m
\Hom_{A_{q_m}}(C_m,J_m) \rightarrow 0\ .
$$
It is a well known fact that the surjectivity of the second map is
equivalent to the vanishing of
${\mathop{\lim}\limits_{\longleftarrow}}^{(1)}$ and hence the
acyclicity for $\mathop{\lim}\limits_{\longleftarrow}$ of the
projective system $(\Hom_{A_{q_m}}(C_m,J_m))_m$.

\medskip

{\bf Lemma 8.2:} {\it For any $C = (C_n)_n \in Coh_{(A,(q_n))}$
and $S = (S_n)_n \in Sh_{(A,(q_n))}$ we have}
$$
\matrix{ \Hom_A(\Gamma(C),\Gamma(S)) & = &
\mathop{\lim\limits_{\longleftarrow}}\limits_{n}
 \Hom_A(\Gamma(C),S_n)\hfill\cr\cr
 & = & \mathop{\lim\limits_{\longleftarrow}}\limits_{n}
 \Hom_{A_{q_n}}(A_{q_n} \otimes_A \Gamma(C),S_n)\hfill\cr\cr
 & = & \Hom_{Sh_{(A,(q_n))}}(C,S)\ .\hfill }
$$

Proof: This is obvious from $C_n = A_{q_n} \otimes_A \Gamma(C)$.

\medskip

{\bf Lemma 8.3:} {\it i. For any two coherent sheaves $C =
(C_n)_n$ and $\widetilde{C} = (\widetilde{C}_n)_n$ we have
$$
\Ext^{\ast}_{Sh_{(A,(q_n))}}(C,\widetilde{C})
\mathop{\rightarrow}\limits^{\cong}
\Ext^{\ast}_A(\Gamma(C),\Gamma(\widetilde{C}))\ ;
$$

ii. for any $C = (C_n)_n \in Coh_{(A,(q_n))}$ and $S = (S_n)_n \in
Sh_{(A,(q_n))}$ we have the short exact sequence}
$$
0 \rightarrow
{\mathop{\lim\limits_{\longleftarrow}}\limits_{n}}^{(1)}
\Ext^{\ast -1}_{A_{q_n}}(C_n,S_n) \rightarrow
\Ext^{\ast}_{Sh_{(A,(q_n))}}(C,S) \rightarrow
\mathop{\lim\limits_{\longleftarrow}}\limits_{n}
\Ext^{\ast}_{A_{q_n}}(C_n,S_n) \rightarrow 0\ .
$$

Proof: Let $\widetilde{C} \longrightarrow J^{\cdot}$ be an
injective resolution in $Sh_{(A,(q_n))}$. Since coherent sheaves
are acyclic for the global section functor (Thm. B) the sequence
$$
\Gamma(\widetilde{C}) \longrightarrow \Gamma(J^{\cdot})
$$
of global sections is exact. As a consequence of Lemma 8.1.iii it
then is an injective resolution of the $A$-module
$\Gamma(\widetilde{C})$. Using Lemma 8.2 we obtain
$$
\matrix{ \Ext^{\ast}_{Sh_{(A,(q_n))}}(C,\widetilde{C}) & = &
h^{\ast}(\Hom_{Sh_{(A,(q_n))}}(C,J^{\ast})) =
h^{\ast}(\Hom_A(\Gamma(C),\Gamma(J^{\ast})))\cr & = &
\Ext^{\ast}_A(\Gamma(C),\Gamma(\widetilde{C}))\ .\hfill }
$$
This establishes the first assertion. For the second one we note
that, by Lemma 8.1.iv, the last identity in Lemma 8.2 gives rise
to a spectral sequence
$$
E^{l,m}_2 =
{\mathop{\lim\limits_{\longleftarrow}}\limits_{n}}^{(l)}
\Ext^m_{A_{q_n}}(C_n,S_n) \Rightarrow
\Ext^{l+m}_{Sh_{(A,(q_n))}}(C,S)\ .
$$
Since the higher derived functors of the projective limit vanish
for $l \geq 2$ this spectral sequence degenerates into the
asserted short exact sequence.

\medskip

Our dimension formalism will be based on the properties of the
functors
$$
\matrix{ \Ext^l_A(.,A)\ :\ {\rm Mod}(A) & \longrightarrow & {\rm
Mod}^r(A) \hfill\cr \hfill M & \longmapsto & \Ext^l_A(M,A) }
$$
for $l \geq 0$.

\medskip

{\bf Lemma 8.4:} {\it For any coadmissible $A$-module $M$ and any
integer $l \geq 0$ the $A$-module $\Ext^l_A(M,A)$ is coadmissible
with
$$
\Ext^l_A(M,A) \otimes_A A_{q_n} = \Ext^l_{A_{q_n}}(A_{q_n}
\otimes_A M,A_{q_n})
$$
for any $n \in \Ndss$.}

Proof: We first claim that
$$
\Ext^l_{A_{q_{n+1}}}(A_{q_{n+1}} \otimes_A M,A_{q_{n+1}})
\otimes_{A_{q_{n+1}}} A_{q_n} = \Ext^l_{A_{q_n}}(A_{q_n} \otimes_A
M,A_{q_n})\ .
$$
But using a resolution of the $A_{q_{n+1}}$-module $A_{q_{n+1}}
\otimes_A M$ by finitely generated free $A_{q_{n+1}}$-modules
together with the flatness of $A_{q_n}$ over $A_{q_{n+1}}$ this
reduces to the obvious equality
$$
\Hom_{A_{q_{n+1}}}(A_{q_{n+1}},A_{q_{n+1}}) \otimes_{A_{q_{n+1}}}
A_{q_n} = \Hom_{A_{q_n}}(A_{q_n},A_{q_n})\ .
$$
It follows that $(\Ext^l_{A_{q_n}}(A_{q_n} \otimes_A
M,A_{q_n}))_n$, for each $l \geq 0$, is a coherent sheaf.
Moreover, Lemma 8.3 and Theorem B imply that
$$
\Ext^l_A(M,A) \mathop{\longrightarrow}\limits^{\cong}
\Gamma((\Ext^l_{A_{q_n}}(A_{q_n}\otimes_A M,A_{q_n})_n)\ .
$$
This isomorphism is $A$-linear as can be seen by computing both
sides using a projective resolution of the $A$-module $M$ (and
observing that $A_{q_n}$ is flat over $A$). It remains to apply
Cor. 3.1.

\medskip

At this point we briefly recall the notion of Auslander regularity
(compare [Bjo] or [LVO] Chap. III). Let $R$ be an arbitrary
associative unital ring. The {\it grade} of a (left or right)
$R$-module $N$ is defined by
$$
j_R(N) := {\rm min}\{l \geq 0 : \Ext^l_R(N,R) \neq 0\}\ .
$$
The $R$-module $N$ is called {\it pure} if
$\Ext^l_R(\Ext^l_R(N,R),R) = 0$ for any $l \neq j_R(N)$. Suppose
now that $R$ is left and right noetherian. If $N \neq 0$ is
finitely generated then its grade $j_R(N)$ is bounded above by the
projective dimension of $N$. One says that $N$ satisfies the {\it
Auslander condition} if, for each $l \geq 0$ and any submodule $L
\subseteq \Ext^l_R(N,R)$, we have $j_R(L) \geq l$. The noetherian
ring $R$ is called (left and right) {\it Auslander regular} if
every finitely generated left or right $R$-module has finite
projective dimension and satisfies the Auslander condition.
Suppose henceforth that $R$ is Auslander regular, and let $N$ be a
finitely generated $R$-module of projective dimension $d(N)$. For
any submodule $N_0 \subseteq N$ we have
$$
j_R(N) = {\rm min}(j_R(N_0),j_R(N/N_0))\ .
$$
Most importantly, $N$ carries a natural filtration, called the
{\it dimension filtration}, by submodules
$$
N = \Delta^0(N) \supseteq \Delta^1(N) \supseteq\ldots\supseteq
\Delta^{d(N)+1}(N) = 0\ .
$$
This filtration is characterized by the property that a submodule
$L \subseteq N$ has grade $j_R(L) \geq l$ if and only if $L
\subseteq \Delta^l(N)$. In addition one has:

-- $j_R(N) = {\rm sup}\{l \geq 0 : \Delta^l(N) = N\}$\ ;

-- if $N$ is pure then $N = \Delta^{j_R(N)}(N) \supseteq
\Delta^{j_R(N)+1}(N) = 0$\ ;

-- $\Delta^l(N)/\Delta^{l+1}(N)$ is zero or pure of grade $l$\ .

These properties show that over an Auslander regular ring the
grade can be viewed as a substitute for the codimension of a
finitely generated module.

\smallskip

Going back to our Fr\'echet-Stein algebra $A$ we assume from now
on that
$$
\matrix{\hbox{there is an integer}\ d \geq 0\ \hbox{such that
each}\ A_{q_n}\ \hbox{is Auslander regular}\cr\hbox{of a global
dimension}\ \leq d\ .\hfill} \leqno{({\rm DIM})}
$$
As an immediate consequence of Lemma 8.4 we obtain for any
coadmissible $A$-module $M$ that:

-- $j_{A_{q_{n+1}}}(A_{q_{n+1}} \otimes_A M) \leq
j_{A_{q_n}}(A_{q_n} \otimes_A M)$ for any $n \in \Ndss$ ;

-- $j_A(M) = \mathop{\rm min}\limits_n\, j_{A_{q_n}}(A_{q_n}
\otimes_A M)$\ ;

-- $j_A(M) \leq d$ if $M \neq 0$ .

In view of Lemma 3.6 it further follows that
$$
j_A(M) = {\rm min}(j_A(N),j_A(M/N))\ \ \hbox{for any coadmissible
submodule}\ N \subseteq M\ .
$$
Because of Lemma 8.4 we may apply the latter to any coadmissible
submodule $L \subseteq \Ext^l_A(M,A)$ and obtain
$$
\matrix{ j_A(L) & \geq & j_A(\Ext^l_A(M,A)) = \mathop{\rm
min}\limits_n\,
 j_{A_{q_n}}(\Ext^l_A(M,A) \otimes_A A_{q_n}) \cr\cr
 & = & \mathop{\rm min}\limits_n\, j_{A_{q_n}} (\Ext^l_{A_{q_n}}(A_{q_n}
 \otimes_A M,A_{q_n})) \hfill\cr\cr & \geq & l\ . \hfill }
$$
This means that in an obvious categorical sense any object in the
category $\Cscr_A$ satisfies the Auslander condition.

\medskip

{\bf Remark 8.5:} {\it For any two $M_1,M_2 \in \Cscr_A$ we have
$$
\Ext^l_A(M_1,M_2) = 0\ \ \ \hbox{for}\ \ l > d+1\ ;
$$
if $M_2$ is finitely generated then}
$$
\Ext^{d+1}_A(M_1,M_2) = 0\ .
$$

Proof: The first vanishing statement is immediate from Lemma 8.3.
For the second one we suppose that $M_2$ is finitely generated.
Writing $M_2$ as a quotient of a finitely generated free module
and using the long exact Ext-sequence together with additivity
reduces us to the case $M_2 = A$ which is a consequence of Lemma
8.4.

\medskip

The possibility of extending the dimension filtration to
coadmissible $A$-modules relies on the following fact.

\medskip

{\bf Lemma 8.6:} {\it Let $R_0 \longrightarrow R_1$ be a unital
homomorphism between two Auslander regular noetherian rings such
that $R_1$ is flat as a left as well as a right $R_0$-module; for
any finitely generated (left or right) $R_0$-module $N$ we have
$$
\Delta^l(R_1 \otimes_{R_0} N) = R_1 \otimes_{R_0} \Delta^l(N)\ \ \
\hbox{for any}\ l \geq 0\ .
$$ }

Proof: There is a convergent spectral sequence (see [Bjo] or
[LVO]) with $E_2$-term $E_2^{l,m} :=
\Ext^l_{R_0}(\Ext^{-m}_{R_0}(N,R_0),R_0)$ and abutment $E^{l+m} :=
N$ for $l+m = 0$, resp. $:= 0$ for $l+m \neq 0$. The filtration on
$N$ induced by this spectral sequence coincides with the dimension
filtration. By using projective resolutions by finitely generated
projective modules it follows from our flatness assumption that
$$
R_1 \otimes_{R_0} \Ext^l_{R_0}(\Ext^{-m}_{R_0}(N,R_0),R_0) =
\Ext^l_{R_1}(\Ext^{-m}_{R_1}(R_1 \otimes_{R_0} N,R_1),R_1)\ .
$$
In fact the base extension $R_1 \otimes_{R_0} .$ transforms the
whole spectral sequence and hence the dimension filtration for the
$R_0$-module $N$ into the corresponding spectral sequence and
dimension filtration for the $R_1$-module $R_1 \otimes_{R_0} N$.

\medskip

Lemma 8.6 implies that any coadmissible $A$-module $M$ carries a
{\it dimension filtration}
$$
M = \Delta^0(M) \supseteq \Delta^1(M) \supseteq\ldots\supseteq
\Delta^{d+1}(M) = 0\ .
$$
by coadmissible submodules $\Delta^l(M)$ such that
$$
A_{q_n} \otimes_A \Delta^l(M) = \Delta^l(A_{q_n} \otimes_A M)\ \ \
\hbox{for any}\ n \in \Ndss\ .
$$

\medskip

{\bf Proposition 8.7:} {\it Let $M$ be a coadmissible $A$-module;
we then have:

i. A coadmissible submodule $N \subseteq M$ has grade $j_A(N) \geq
l$ if and only if $N \subseteq \Delta^l(M)$;

ii. $j_A(M) = {\rm sup}\{l \geq 0 : \Delta^l(M) = M\}$\ ;

iii. if $M$ is pure then $M = \Delta^{j_A(M)}(M) \supseteq
\Delta^{j_A(M)+1}(M) = 0$\ ;

iv. $\Delta^l(M)/\Delta^{l+1}(M)$ is zero or pure of grade $l$\ .}

Proof: i. If $N \subseteq \Delta^l(M)$ then $j_A(N) \geq
j_A(\Delta^l(M)) = {\rm min} j_{A_{q_n}}(\Delta^l(A_{q_n}
\otimes_A M)) \geq l$. On the other hand, if $j_A(N) \geq l$ then
$j_{A_{q_n}}(A_{q_n} \otimes_A N) \geq l$ for any $n$, hence
$A_{q_n} \otimes_A N \subseteq \Delta^l(A_{q_n} \otimes_A M)$ for
any $n$, and therefore $N \subseteq \Delta^l(M)$. ii. This is a
formal consequence of the first assertion. iii. We have $j_A(M) =
j_{A_{q_n}}(A_{q_n} \otimes_A M)$ for $n \in \Ndss$ big enough.
Lemma 8.4 then implies, for these $n$, that $A_{q_n} \otimes_A M$
is pure and hence that $A_{q_n} \otimes_A M =
\Delta^{j_A(M)}(A_{q_n} \otimes_A M) \supseteq
\Delta^{j_A(M)+1}(A_{q_n} \otimes_A M) = 0$. iv. Suppose that
$\Delta^l(M)/\Delta^{l+1}(M)$ is nonzero. According to Lemma 8.4
it suffices to show that $A_{q_n} \otimes_A
(\Delta^l(M)/\Delta^{l+1}(M))$ is pure of grade $l$ for any $n \in
\Ndss$ big enough. Since $A_{q_n}$ is flat over $A$ (Remark 3.2)
we have $A_{q_n} \otimes_A (\Delta^l(M)/\Delta^{l+1}(M)) =
(A_{q_n} \otimes_A \Delta^l(M))/(A_{q_n} \otimes_A
\Delta^{l+1}(M)) = \Delta^l(A_{q_n} \otimes_A
M)/\Delta^{l+1}(A_{q_n} \otimes_A M)$. By assumption the left hand
side is nonzero for big $n$. If so the right hand side is pure of
grade $l$.

\medskip

All these properties confirm our claim that under the assumption
(DIM) the grade can be used as a well behaved codimension function
on the category $\Cscr_A$.

Before we show that all this applies to the algebra $D(G,K)$ for a
compact locally $\Qdss_p$-analytic group $G$ we recall, for the
convenience of the reader, the following facts from general ring
theory.

\medskip

{\bf Lemma 8.8:} {\it Let $R_0 \longrightarrow R_1$ be a unital
homomorphism of (left or right) noetherian rings; suppose that
there are units $b_1=1,b_2,\ldots,b_m \in (R_1)^{\times}$ which
form a basis of $R_1$ as an $R_0$-module and which
satisfy:\hfill\break
 -- $b_iR_0 = R_0b_i$ for any $1 \leq i \leq m$,\hfill\break
 -- for any $1 \leq i,j \leq m$ there is a $1 \leq k \leq m$ such
 that $b_ib_j \in b_kR_0$, and\hfill\break
 -- for any $1 \leq i \leq m$ there is a $1 \leq l \leq m$ such
 that $b_i^{-1} \in b_lR_0$.\hfill\break
We then have
$$
\Ext^{\ast}_{R_1}(N_1,R_1\otimes_{R_0}N_0),\ resp.\
\Ext^{\ast}_{R_1}(N_1,N_0\otimes_{R_0}R_1) \cong
\Ext^{\ast}_{R_0}(N_1,N_0)\leqno{(+)}
$$
for any pair of finitely generated left, resp. right,
$R_i$-modules $N_i$ ($i = 0,1$). If $N$ is any finitely generated
(left or right) $R_1$-module then:\hfill\break
 i. $\Ext^{\ast}_{R_1}(N,R_1) \cong \Ext^{\ast}_{R_0}(N,R_0)$ and in
particular $j_{R_1}(N) = j_{R_0}(N)$;\hfill\break
 ii. the projective dimensions of $N$ as an $R_1$-module and
as an $R_0$-module coincide; in particular, $N$ is projective as
an $R_1$ -module if and only if it is projective as an
$R_0$-module.\hfill\break
 Finally, the rings $R_0$ and $R_1$ have the same global dimension.}

Proof: The case of right modules being analogous we only consider
the case of left modules in the following. Let $\ell : R_1
\longrightarrow R_0$ denote the projection onto the first summand
in the decomposition
$$
R_1 = \bigoplus_i b_iR_0 = \bigoplus_i R_0b_i\ .
$$
It induces a map between the two sides in $(+)$. By using a
projective resolution of the $R_1$-module $N_1$ by finitely
generated free $R_1$-modules (which as a consequence of our
assumptions also is such a resolution for $N_1$ as an
$R_0$-module) it suffices to consider the case $\ast = 0$ and $N_1
= R_1$. By using a presentation of $N_0$ by finitely generated
free $R_0$-modules we are further reduced to show that the map
$$
\matrix{ \Hom_{R_1}(R_1,R_1) & \longrightarrow &
\Hom_{R_0}(R_1,R_0)\cr \hfill \Phi & \longmapsto & \ell\circ\Phi
\hfill }
$$
is bijective. We have
$$
\Phi(a) = b_1(\ell\circ\Phi)(b_1^{-1}a) + \ldots +
b_m(\ell\circ\Phi)(b_m^{-1}a) \qquad \hbox{for any}\ a \in R_1
$$
which proves the injectivity of the map. To establish surjectivity
we suppose given a left $R_0$-module map $\Psi : R_1
\longrightarrow R_0$. It is straightforward to check that
$$
\Phi(a) := b_1\Psi(b_1^{-1}a) + \ldots + b_m\Psi(b_m^{-1}a)
$$
in fact is a preimage of $\Psi$. This proves the isomorphism
$(+)$.

The statement i. is $(+)$ for $N_0 = R_0$. In the statement ii. we
first consider the particular case. If $N$ is projective as an
$R_1$-module then obviously also as an $R_0$-module since $R_1$ is
finitely generated free over $R_0$. Suppose therefore that $N$ is
projective over $R_0$. We have to show that $\Ext^1_{R_1}(N,X) =
0$ for all $R_1$-modules $X$. Since $R_1$ is noetherian and $N$ is
finitely generated the functor $\Ext^1_{R_1}(N,.)$ commutes with
filtered inductive limits. Hence it suffices to consider finitely
generated modules $X$. From $(+)$ we then have that
$\Ext^1_{R_1}(N,R_1\otimes_{R_0} X) = 0$. But $X$ as an
$R_1$-module is a direct summand of $R_1\otimes_{R_0} X$ since
$R_0$ is a direct $R_0$-module summand of $R_1$. The general
assertion about the equality of projective dimensions as well as
the asserted equality of global dimensions is easily deduced from
this using appropriate projective resolutions (compare [MCR] Thm.
7.5.6).

\medskip

{\bf Theorem 8.9:} {\it For any compact locally $\Qdss_p$-analytic
group $G$ the Fr\'echet-Stein algebra $D(G,K)$ satisfies (DIM)
with $d := {\rm dim}(G)$.}

Proof: In a first step we assume in addition that $G$ is a uniform
pro-$p$-group. The structure of $D(G,K)$ as a two sided
Fr\'echet-Stein algebra is given, according to section 4, by the
Banach algebras $D_r(G,K)$ where $1/p < r < 1, r \in p^{\Qdss}$.
We have to show that each $D_r(G,K)$ is Auslander regular of
global dimension $\leq d$. The filtration on $D_r(G,K)$ derived
from the submultiplicative norm $\|\ \|_r$ is quasi-integral. From
Thm. 4.5.i we know that the associated graded ring $gr^\cdot_r
D_r(G,K)$ is isomorphic to a polynomial ring
$k[X_0,X_0^{-1},X_1,\ldots,X_d]$ over the residue class field $k$
of $K$. This is a regular commutative noetherian ring of global
dimension $d+1$. By [LVO] Cor. I.7.2.2, Thm. III.2.2.5, and
III.2.4.3 it follows that $D_r(G,K)$ is Auslander regular of
global dimension $\leq d+1$. To obtain the stronger bound $\leq d$
we have to reexamine part of this reasoning more closely. By [MCR]
7.2.6 we may replace $K$ by a finite extension. This allows us, by
the computation in the proof of Lemma 4.8, to assume that
$gr^\cdot_r F^0_rD_r(G,K)$ is isomorphic to a polynomial ring
$k[Y_0,\ldots,Y_d]$. Hence the global dimension of $F^0_rD_r(G,K)$
also is $\leq d+1$. On the other hand since $F^{0+}_rD_r(G,K)$ is
contained in the Jacobson radical of $F^0_rD_r(G,K)$ (see [LVO]
Lemma I.3.5.5(2)) multiplication by $p$ is zero on any simple
module over $F^0_rD_r(G,K)$. It therefore follows from [MCR]
7.4.3/4 that the global dimension of $D_r(G,K) = \Qdss_p \otimes
F^0_rD_r(G,K)$ is $\leq d$.

For general $G$ we choose an open normal subgroup $H \subseteq G$
which is a uniform pro-$p$-group. In view of the construction of
the algebras $D_r(G,K)$ in the proof of Thm. 5.1 the ring
extensions $D_r(H,K) \subseteq D_r(G,K)$ satisfy the assumptions
of Lemma 8.8. It follows that with $D_r(H,K)$ also $D_r(G,K)$ is
Auslander regular of global dimension $\leq d$.

\medskip

{\bf Remark 8.10:} {\it Let $H$ be an open subgroup in the compact
locally $\Qdss_p$-analytic group $G$; for any coadmissible
$D(G,K)$-module $M$ we have}
$$
j_{D(G,K)}(M) = j_{D(H,K)}(M)\ .
$$

Proof: This follows from Lemma 8.8 applied to the ring extensions
$D_r(H_1,K) \subseteq D_r(H,K)$ and $D_r(H_1,K) \subseteq
D_r(G,K)$ provided we find a uniform pro-$p$-group $H_1$ contained
in $H$ which is open normal in $G$. But it is a straightforward
consequence of [DDMS] Thm. 4.2 and the discussion before Thm. 4.9
that $G$ has a fundamental system of open characteristic subgroups
which are uniform pro-$p$.

\medskip

We finally consider an arbitrary locally $\Qdss_p$-analytic group
$G$ of dimension $d$ and a coadmissible (left) $D(G,K)$-module
$M$. Remark 8.10 implies that, for any two compact open subgroups
$H,H' \subseteq G$, we have
$$
j_{D(H,K)}(M) = j_{D(H',K)}(M)\ .
$$
This fact allows us to unambiguously define the {\it codimension}
of $M$ by
$$
{\rm codim}(M) := j_{D(H,K)}(M)\ .
$$
Similarly, by Remark 8.10 and Prop. 8.7, the dimension filtration
$$
M = \Delta^0(M) \supseteq \Delta^1(M) \supseteq\ldots\supseteq
\Delta^{d+1}(M) = 0\ .
$$
of $M$ as a $D(H,K)$-module does not depend on the choice of $H$.

\goodbreak

{\bf Proposition 8.11:} {\it For any coadmissible $D(G,K)$-module
$M$ we have:

i. Each $\Delta^l(M)$ is a coadmissible $D(G,K)$-submodule of $M$;

ii. a coadmissible $D(G,K)$-submodule $N \subseteq M$ has
codimension $\geq l$ if and only if $N \subseteq \Delta^l(M)$;

iii. ${\rm codim}(M) = {\rm sup}\{l \geq 0 : \Delta^l(M) = M\}$\ ;

iv. all nonzero coadmissible $D(G,K)$-submodules of
$\Delta^l(M)/\Delta^{l+1}(M)$ have codimension $l$\ .}

Proof: Everything except the $D(G,K)$-invariance of $\Delta^l(M)$
is a consequence of Prop. 8.7. Since $D(G,K)$ as a $D(H,K)$-module
is generated by the $\delta_g$ for $g \in G$ it suffices to show
that $g(\Delta^l(M)) \subseteq \Delta^l(M)$. The map $\Delta^l(M)
\mathop{\longrightarrow}\limits^{g\cdot} g(\Delta^l(M))$ is a
module isomorphism relative to the algebra isomorphism $\delta_g .
\delta_{g^{-1}} : D(H,K) \mathop{\longrightarrow}\limits^{\cong}
\break D(gHg^{-1},K)$. It follows that $l \leq
j_{D(H,K)}(\Delta^l(M)) = j_{D(gHg^{-1},K)}(g(\Delta^l(M))$. Since
the dimension filtration is independent of the choice of the
compact open subgroup we obtain $g(\Delta^l(M)) \subseteq
\Delta^l(M)$.

\medskip

{\bf Definition:} {\it A coadmissible $D(G,K)$-module $M$ with
$codim(M) \geq dim(G)$ is called zero-dimensional.}

\medskip

A nonzero zero-dimensional $D(G,K)$-module $M$ satisfies ${\rm
codim}(M) = {\rm dim}(G)$.

\medskip

{\bf Theorem 8.12:} {\it Let $G$ be a locally $\Qdss_p$-analytic
group and let $M$ be a coadmissible $D(G,K)$-module; if the action
of the universal enveloping algebra $U(\gfr)$ on $M$ is locally
finite, i.e., if $U(\gfr)x$, for any $x \in M$, is a finite
dimensional $\Qdss_p$-vector space then $M$ is zero-dimensional.}

\medskip

{\bf Corollary 8.13:} {\it If $V$ is an admissible smooth
$G$-representation then $V'_b$ is a zero-dimensional
$D(G,K)$-module.}

\medskip

The statement of Thm. 8.12 only depends on an arbitrarily chosen
compact open subgroup of $G$. We therefore let $H$ in the
following be a compact locally $\Qdss_p$-analytic group which is a
uniform pro-$p$-group. In this case Thm. 8.12 will be a
consequence of the following more general, but technical
criterion. Recalling some notation from section 4 we fix an
ordered basis $(h_1,\ldots,h_d)$ of $G$ and let $\psi : \Zdss_p^d
\mathop{\longrightarrow}\limits^{\sim} G$ given by
$\psi((x_1,\ldots x_d)) = h_1^{x_1}\cdot\ldots\cdot h_d^{x_d}$
denote the corresponding global chart. Any distribution $\lambda
\in D(H,K)$ has a unique convergent expansion in the elements $b_1
= h_1 - 1,\ldots,b_d = h_d - 1$. As always we let $r$ vary over
the set $\{r \in p^{\Qdss} : 1/p < r < 1\}$.

\medskip

{\bf Proposition 8.14:} {\it Suppose that
$\lambda_1,\ldots,\lambda_d \in D(H,K)$ are nonzero elements such
that, for each $1 \leq i \leq d$, the expansion of $\lambda_i$
only involves powers of $b_i$, and let $J \subseteq D(H,K)$ be the
left ideal generated by $\lambda_1,\ldots,\lambda_d$; then the
coadmissible $D(H,K)$-module $D(H,K)/J$ is zero-dimensional.}

Proof: Since $J$ is finitely generated the module $D(H,K)/J$ is
coadmissible and the corresponding coherent sheaf is given by
$D_r(H,K)/D_r(H,K)J$. To simplify notation we use in the following
the abbreviation $D_r := D_r(H,K)$. So we have to show that
$j_{D_r}(D_r/D_rJ) \geq d$. We equip $D_rJ$ and $D_r/D_rJ$ with
the filtrations induced by the filtration $F_r^{\cdot}D_r$. Hence
$gr_r^\cdot(D_r/D_rJ) = (gr_r^\cdot D_r)/(gr_r^\cdot D_rJ)$.
Recall from Thm. 4.5.i that $gr_r^\cdot D_r = gr^\cdot
K[\sigma_r(b_1),\ldots,\sigma_r(b_d)]$ is a polynomial ring in $d$
variables over $gr^\cdot K$. According to [LVO] Thm. III.2.5.2 we
have $j_{D_r}(D_r/D_rJ) = j_{gr_r^\cdot D_r}(gr_r^\cdot
(D_r/D_rJ))$. We therefore are reduced to showing that
$j_{gr_r^\cdot D_r}(gr_r^\cdot D_r/gr_r^\cdot D_rJ) \geq d$. The
ring $ R := gr_r^\cdot D_r$ is a commutative (Auslander) regular
and catenary noetherian domain of Krull dimension $d+1$. For any
finitely generated module $N$ over such a ring $R$ standard
commutative algebra (compare [BH] Cor. 3.5.11) implies the formula
$$
j_R(N) = d+1 - {\rm Krulldim}(R/{\rm ann}(N))
$$
where ${\rm ann}(N)$ denotes the annihilator ideal of $N$ in $R$.
Hence it remains to see that the factor ring $gr_r^\cdot
D_r/gr_r^\cdot D_rJ$ has Krull dimension $\leq 1$. By our
assumption on the $\lambda_i$ we have a surjection
$$
gr^\cdot K[\sigma_r(b_1)]/\langle\sigma_r(\lambda_1)\rangle
\mathop{\otimes}\limits_{gr^\cdot K} \ldots
\mathop{\otimes}\limits_{gr^\cdot K} gr^\cdot
K[\sigma_r(b_d)]/\langle\sigma_r(\lambda_d)\rangle
\dlongrightarrow gr_r^\cdot D_r/gr_r^\cdot D_rJ\ .
$$
Each $gr^\cdot K[\sigma_r(b_i)]/\langle\sigma_r(\lambda_i)\rangle$
is finitely generated as a $gr^\cdot K$-module. It follows that
$gr_r^\cdot D_r/gr_r^\cdot D_rJ$ is finite over the one
dimensional ring $gr^\cdot K$ and therefore has Krull dimension
$\leq 1$ (compare [B-CA] V\S2.1).

\medskip

$\underline{\rm Proof\ of\ Thm.\, 8.12}:$ According to Prop. 8.11
we have to show that $M \subseteq \Delta^d(M)$. In fact, it
suffices to show that, for any $x \in M$, the submodule $N :=
D(H,K)x$ of $M$ generated by $x$ (which is coadmissible by Cor.
3.4.iv) is zero-dimensional. We write $N = D(H,K)/J_1$ for some
left ideal $J_1 \subseteq D(H,K)$. By assumption $J_1$ contains an
ideal of finite codimension $J_0$ of $U(\gfr)$. The Lie algebra
$\gfr$ has the basis $\partial_1,\ldots,\partial_d$ where
$\partial_i := \psi_{\ast}((\partial/\partial x_i)_{|x_i=0}) =
{\rm log}(1+b_i)$. Since $J_0$ is of finite codimension in
$U(\gfr)$ it contains, for any $1 \leq i \leq d$, a nonzero
polynomial $P_i(\partial_i) \in K[\partial_i]$. The left ideal $J
\subseteq J_1 \subseteq D(H,K)$ generated by
$P_1(\partial_1),\ldots,P_d(\partial_d)$ satisfies the assumptions
of Prop. 8.14. Hence $D(H,K)/J$ and a fortiori its quotient $N$
are zero-dimensional.

\medskip

The subsequent criterion is obtained by a very similar argument.

\medskip

{\bf Theorem 8.15:} {\it Let $G$ be a compact locally
$\Qdss_p$-analytic group and let $M$ be a coadmissible
$D(G,K)$-module; if the coherent sheaf corresponding to $M$ is a
sheaf of finite dimensional $K$-vector spaces then $M$ is
zero-dimensional.}

Proof: Using Remark 8.10 and the notations introduced before Prop.
8.14 the assertion reduces to the claim that $j_{D_r(H,K)}(M) \geq
d$ for any $D_r(H,K)$-module $M$ which is finite dimensional as a
$K$-vector space. We equip $M$ with a filtration $F^\cdot M$ which
is good with respect to $F_r^\cdot D_r(H,K)$ and obtain from [LVO]
Thm. III.2.5.2 that $j_{D_r(H,K)}(M) = j_{gr_r^\cdot
D_r(H,K)}(gr^\cdot M)$. The assumption that $M$ is finite
dimensional over $K$ easily implies that $gr^\cdot M$ is finitely
generated over $gr^\cdot K$. From this point on the rest of the
argument is exactly the same as in the proof of Prop. 8.14.

\medskip

We strongly believe that the results of this section extend to
locally $L$-analytic groups $G$ over any finite extension $L$ of
$\Qdss_p$. But to establish the Auslander regularity (or possibly
a weaker but sufficient Auslander-Gorenstein property - compare
[Bjo]) of the corresponding Banach algebras in this situation
seems to require new ideas.

\bigskip

{\bf References}

\parindent=23truept

\ref{[Ban]} B\u{a}nic\u{a} C.: Une caracterisation de la dimension
d'un faisceau analitique coherent. Compositio math. 25, 101-108
(1972)

\ref{[Bjo]} Bj{\"o}rk J.-E.: Filtered Noetherian Rings. In
Noetherian rings and their applications, Math. Survey Monographs
24, pp. 59-97, AMS 1987

\ref{[BGR]} Bosch S., G\"untzer U., Remmert R: Non-Archimedean Analysis.
Ber-lin-Heidelberg-New York: Springer 1984

\ref{[B-CA]} Bourbaki N.: Commutative Algebra. Paris: Hermann 1972

\ref{[B-GT]} Bourbaki N.: General Topology, Chap. 1-4. Berlin-Heidelberg-New York:
Springer 1989

\ref{[B-TVS]} Bourbaki N.: Topological Vector Spaces.
Berlin-Heidelberg-New York: Springer 1987

\ref{[BH]} Bruns W., Herzog J.: Cohen-Macaulay rings. Cambridge
Univ. Press 1993

\ref{[DDMS]} Dixon J.D., du Sautoy M.P.F., Mann A., Segal D.:
Analytic Pro-$p$-Groups. Cambridge Univ. Press 1999

\ref{[Fea]} F\'eaux de Lacroix C. T.: Einige Resultate \"uber die topologischen
Dar- stellungen $p$-adischer Liegruppen auf unendlich
dimensionalen Vektor- r\"aumen \"uber einem $p$-adischen K\"orper.
Thesis, K\"oln 1997, Schriftenreihe Math. Inst. Univ. M\"unster,
3. Serie, Heft 23, pp. 1-111 (1999)

\ref{[For]} Forster O.: Zur Theorie der Steinschen Algebren und
Moduln. Math. Z. 97, 376-405 (1967)

\ref{[EGA]} Grothendieck A., Dieudonn\'e J.: \'El\'ements de
g\'eom\'etrie alg\'ebrique, Chap. III. Publ. Math. IHES 11 (1961)

\ref{[Laz]} Lazard M.: Groupes analytiques $p$-adique. Publ. Math.
IHES 26, 389-603 (1965)

\ref{[LVO]} Li Huishi, van Oystaeyen F.: Zariskian Filtrations.
Dordrecht: Kluwer 1996

\ref{[MCR]} McConnell J.C., Robson J.C.: Noncommutative noetherian
rings. Chi- chester: Wiley 1987

\ref{[Pop]} Popescu N.: Abelian Categories with Applications to
Rings and Mo\-dules. London-New York: Academic Press 1973

\ref{[NFA]} Schneider P.: Nonarchimedean Functional Analysis. Berlin-Heidelberg-New York:
Springer 2001

\ref{[ST1]} Schneider P., Teitelbaum J.:  $U(\gfr)$-finite locally analytic
representations. Representation Theory 5, 111-128 (2001)

\ref{[ST2]} Schneider P., Teitelbaum J.: Locally analytic
distributions and $p$-adic representation theory, with
applications to $GL_2$. J. AMS 15, 443-468 (2002)

\ref{[ST3]} Schneider P., Teitelbaum J.: Banach space
representations and Iwasawa theory. Israel J. Math. 127, 359-380
(2002)

\ref{[ST4]} Schneider P., Teitelbaum J.: $p$-adic Fourier theory.
Documenta Math. 6, 447-481 (2001)

\ref{[ST5]} Schneider P., Teitelbaum J.: $p$-adic boundary values.
In Cohomologies $p$-adiques et applications arithm\'etiques (I)
(Eds. Berthelot/Fontaine/ Illusie/Kato/Rapoport), Ast\'erisque
278, 51-125 (2002)

\bigskip

\parindent=0pt

Peter Schneider\hfill\break Mathematisches Institut\hfill\break
Westf\"alische Wilhelms-Universit\"at M\"unster\hfill\break
Einsteinstr. 62\hfill\break D-48149 M\"unster, Germany\hfill\break
pschnei@math.uni-muenster.de\hfill\break
http://www.uni-muenster.de/math/u/schneider\hfill

\noindent
Jeremy Teitelbaum\hfill\break Department of Mathematics, Statistics,
and Computer Science (M/C 249)\hfill\break University of Illinois at
Chicago\hfill\break 851 S. Morgan St.\hfill\break Chicago, IL 60607,
USA\hfill\break jeremy@uic.edu\hfill\break
http://raphael.math.uic.edu/$\sim$jeremy\hfill

\end